# On the Virtual Element Method for Three-Dimensional Elasticity Problems on Arbitrary Polyhedral Meshes

Arun L. Gain†, Cameron Talischi†, Glaucio H. Paulino*

*Department of Civil & Environmental Engineering, University of Illinois at Urbana-Champaign, USA*

†Equal contribution authors,

* *Corresponding author:* paulino@illinois.edu

*November 4, 2013*

**Abstract**

We explore the recently-proposed Virtual Element Method (VEM) for numerical solution of boundary value problems on arbitrary polyhedral meshes. More specifically, we focus on the elasticity equations in three-dimensions and elaborate upon the key concepts underlying the first-order VEM. While the point of departure is a conforming Galerkin framework, the distinguishing feature of VEM is that it does not require an explicit computation of the trial and test spaces, thereby circumventing a barrier to standard finite element discretizations on arbitrary grids. At the heart of the method is a particular kinematic decomposition of element deformation states which, in turn, leads to a corresponding decomposition of strain energy. By capturing the energy of linear deformations exactly, one can guarantee satisfaction of the engineering patch test and optimal convergence of numerical solutions. The decomposition itself is enabled by local projection maps that appropriately extract the rigid body motion and constant strain components of the deformation. As we show, computing these projection maps and subsequently the local stiffness matrices, in practice, reduces to the computation of purely geometric quantities. In addition to discussing aspects of implementation of the method, we present several numerical studies in order to verify convergence of the VEM and evaluate its performance for various types of meshes.

*Keywords*: Virtual Element Method, Mimetic Finite Differences, polyhedral meshes, polytopes, Voronoi tessellations

## 1. Introduction

The development of discretization methods for solving three-dimensional boundary value problems on general polyhedral meshes has recently received considerable attention in the numerical analysis literature. One driving force behind this trend is the difficulty associated with mesh generation for complex or evolving domains, for which the use of arbitrarily-shaped elements can provide much needed flexibility [22]. For example, a simple embedding strategy consisting of carving out the problem domain out of a structure background grid, produces polyhedral elements at the boundary [31]. Mesh refinement and coarsening in adaptive schemes can also be handled with greater ease if the analysis allows for the presence of elements with general geometries [30]. In addition to advantages in mesh generation and adaptation, polyhedral discretizations can deliver improved performance in some applications. For example, as discussed in [21], polyhedral meshes can achieve the same level of accuracy in flow simulations compared to their simplicial counterparts but with far fewer number of cells and unknowns.

With regards to the type of discretization method, finite volume methods based on polyhedral cells have reached a level of maturity in fluid dynamic simulations, as evidenced by their availability and use in commercial software [1, 2]. Mimetic finite difference (MFD) methods, capable of handling general three-dimensional meshes, are also the subject of active research and have been successfully applied to diffusion, elasticity, and fluid flow problems (see, for example, [16, 17, 18, 8, 5]). The extension of finite element methods in this arena, however, has been relatively slow, despite the availability of special interpolation functions in the literature. This is, in part, due to the fact that these interpolants are subject to restrictions on the topology of admissible elements (e.g., convexity, maximum valence count) and can be sensitive to geometric degeneracies. More importantly, calculating these functions and their gradients are often prohibitively expensive. Numerical evaluation of weak form integrals, with sufficient accuracy, poses yet another challenge due to the non-polynomial nature of these functions as well as the arbitrary domain of integration[1] [37]. To mention a few approaches in the literature aiming to overcome these barriers, we point to the work by Rashid and co-workers [31, 30], who have developed elements based on non-conforming polynomial or piecewise polynomial basis functions that are tolerant of degeneracies. More recently, harmonic basis functions have been considered by [32, 14] with particular attention to alleviating the cost of their computation and integration. Other works include constructions based on natural element [29], non-Siboson [24], and mean value coordinates [39].

In this work, we focus on the recently-developed Virtual Element Method (VEM) that addresses some of the above-mentioned challenges facing finite element schemes [10, 11, 15]. As with finite elements, VEM is a Galerkin scheme with an underlying approximation space defined according to a partition (mesh) of the domain. However, it is distinguished from classical finite elements in that it does not require the computation of the interpolation functions in the interior of the elements. One goal of the present work is to break down and elaborate upon the core mathematical concepts underlying VEM within the context of elasticity boundary value problems. The key to the success of the method is a consistent approximation to the elemental strain energy that is exact for the linear deformations without requiring volumetric integration of the basis functions. What enables this approximation is a set of local projection maps that appropriately split up the element deformation into its polynomial and non-polynomial components. In the case of the first-order VEM formulation, where the degrees of freedom are associated with the vertices of the elements, two projection maps, associated with rigid body motion and constant strain deformations, respectively, are used to achieve this kinematic decomposition. As we shall discuss, these projection maps can be beneficial even for finite element schemes when one has access to interpolation functions. We should note that while VEM provides the general recipe for extension to higher-order and higher-continuity polyhedral elements (see, for example, [6]), this significant technology may be hidden in this paper as we will limit the discussion, for the sake of clarity, to the first-order formulation.

[1]By contrast, classical finite elements feature interpolation functions that are either polynomials or images of polynomials and numerical integration is carried out by means of a mapping to a fixed parent domain.

The other task undertaken here is to discuss, in detail, aspects of the implementation of the method for general polyhedral meshes. To this effect, we will derive explicit expressions for the element stiffness matrix and discrete representation of the element projection maps. In addition to two matrices containing special arrangements of coordinates of the element vertices, we encounter two matrices that require calculation of surface integrals of the basis functions over the element boundary. These quantities also reduce to geometric information of the faces (centroids, areas, etc.), if either the approximation spaces are based on interpolants derived by [3] or if a consistent nodal quadrature rule is used. While the connection between MFD method and VEM has been established in the original papers on VEM, the discussion here further elucidates this relationship and illustrates how the Galerkin framework with an underlying approximation space serves as a vehicle for constructing a method that is ultimately geometric in nature.

The remainder of this paper is organized as follows. In Section 2, we define the three-dimensional elasticity model problem and its Galerkin approximation on a polyhedral mesh. Section 3 presents the VEM formulations and its theoretical underpinning. Next, in Section 4, we derive explicit and simplified expressions for the element stiffness matrices and discuss aspects of implementation of the method. We evaluate the performance of VEM in Section 5 via several numerical studies and conclude the work with some remarks in Section 6.

We follow fairly standard notation throughout the paper. As usual, Sobolev spaces $H^k(\Omega)$ consists of functions whose derivatives up to the $k$th order are square-integrable on $\Omega$. The norm on this space, as well as its vector-valued counterpart $H^k(\Omega)^3$, is denoted by $\|\cdot\|_{k,\Omega}$. We denote the symmetric gradient operator by $\boldsymbol{\epsilon}(\cdot) = \left(\nabla \cdot + \nabla^\top \cdot\right)/2$ and skew-symmetric gradient by $\boldsymbol{\omega}(\cdot) = \left(\nabla \cdot - \nabla^\top \cdot\right)/2$. Also, $\boldsymbol{I}_m$ represents the $m \times m$ identity matrix. We shall denote the components of vectors, matrices and tensors in the canonical Euclidean basis with subscripts inside parentheses (e.g. $\boldsymbol{v}_{(i)}$ or $\boldsymbol{\epsilon}_{(ij)}$) in order to make a distinction with indexed quantities. Finally, we use $|\cdot|$ to denote the measure (area or volume) of a set as well as the Euclidian norm of a vector.

## 2. Model problem and discretization

Consider a linear elastic body, with constant stiffness tensor $\boldsymbol{C}$, occupying a smooth bounded domain $\Omega \subseteq \Re^3$ whose boundary $\partial\Omega$ is partitioned into disjoint non-trivial segments $\Gamma_{\boldsymbol{u}}$ and $\Gamma_{\boldsymbol{t}}$. The body is subjected to body forces $\boldsymbol{b}$ in $\Omega$, surface tractions $\boldsymbol{t}$ on $\Gamma_{\boldsymbol{t}}$, and applied displacements $\boldsymbol{g}$ on $\Gamma_{\boldsymbol{u}}$, with all fields assumed to have sufficient regularity (cf. Fig. 1). The resulting deformation $\boldsymbol{u}$ is the unique minimizer of the total potential energy:

$$\boldsymbol{u} = \underset{\boldsymbol{v} \in \mathcal{V}^{\boldsymbol{g}}}{\operatorname{argmin}} \left( \frac{1}{2} a(\boldsymbol{v}, \boldsymbol{v}) - f(\boldsymbol{v}) \right) \tag{1}$$

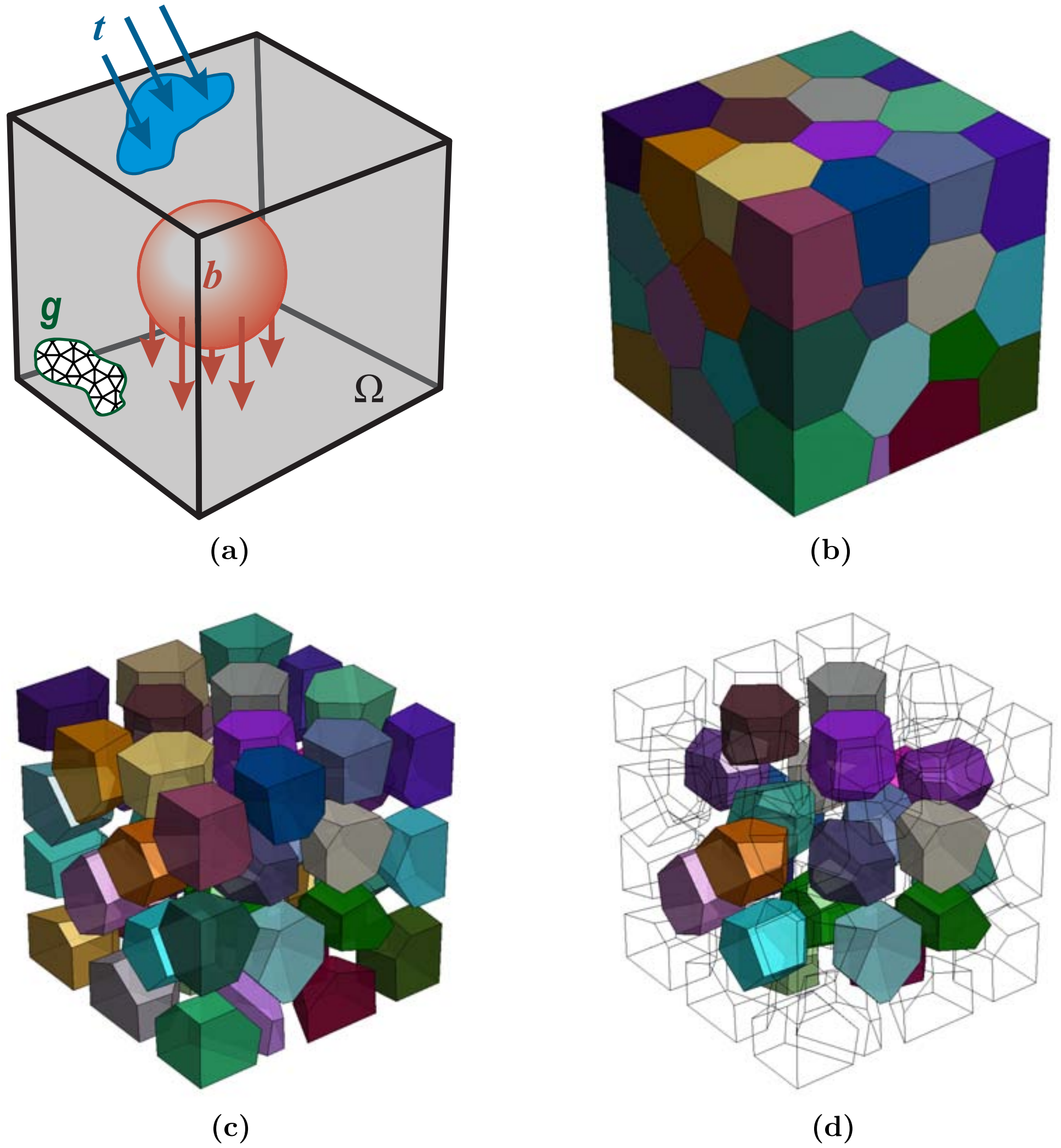

**Figure 1:** (a) Schematic illustrating the elasticity boundary value problem. (b) Partition $\mathcal{T}_h$ of the domain $\Omega$. (c) Split view of the discretized domain. (d) View of a few boundary and internal elements.

Here, $\mathcal{V}^{\boldsymbol{g}} = \{\boldsymbol{v} \in H^1(\Omega)^3 : \boldsymbol{v} = \boldsymbol{g} \text{ on } \Gamma_{\boldsymbol{u}}\}$ is the space of kinematically admissible displacement fields and

$$a(\boldsymbol{u}, \boldsymbol{v}) = \int_\Omega \boldsymbol{\sigma}(\boldsymbol{u}) : \boldsymbol{\epsilon}(\boldsymbol{v}) \mathrm{d}\boldsymbol{x}, \qquad f(\boldsymbol{v}) = \int_\Omega \boldsymbol{b} \cdot \boldsymbol{v} \mathrm{d}\boldsymbol{x} + \int_{\Gamma_{\boldsymbol{t}}} \boldsymbol{t} \cdot \boldsymbol{v} \mathrm{d}\boldsymbol{s} \tag{2}$$

are the energy bilinear form and load linear form, respectively. In the above expression, $\boldsymbol{\sigma}(\boldsymbol{u})$ denotes the stress field associated with $\boldsymbol{u}$, that is,

$$\boldsymbol{\sigma}(\boldsymbol{u}) = \boldsymbol{C}\boldsymbol{\epsilon}(\boldsymbol{u}) \tag{3}$$

Observe that due to the symmetries of $\boldsymbol{C}$, the bilinear form is symmetric in its arguments. We shall assume throughout the paper that the exact solution $\boldsymbol{u}$ is a smooth function (e.g. belongs to $H^2(\Omega)^3$).

### 2.1. Galerkin approximation

Given a sufficiently regular[2] partition $\mathcal{T}_h$ of $\Omega$ into disjoint non-overlapping polyhedra with maximum diameter $h$, we define a conforming discrete space $\mathcal{V}_h$ consisting of continuous displacement fields whose restriction to $E \in \mathcal{T}_h$ belong to the finite-dimensional space $\mathcal{W}(E)$ of smooth functions. Thus,

$$\mathcal{V}_h = \left\{ \boldsymbol{v} \in C^0(\overline{\Omega})^3 : \boldsymbol{v}|_E \in \mathcal{W}(E) \text{ for all } E \in \mathcal{T}_h \right\} \tag{4}$$

Here, the space $\mathcal{W}(E)$ contains the deformation states that can be represented by the element $E$. Note that $C^0$-continuity along with the requirement that $\mathcal{W}(E) \subseteq H^1(E)^3$ implies $\mathcal{V}_h$ is a subspace of $H^1(\Omega)^3$. Moreover, the usual requirement that $\mathcal{W}(E)$ includes linear displacements fields (or equivalently that $E$ can represent rigid body motions and constant states of strains) furnishes the first-order approximation property of $\mathcal{V}_h$, namely that sufficiently smooth displacement fields, including the solution of the continuous problem (1), can be approximated by elements of $\mathcal{V}_h$ with $\mathcal{O}(h)$ errors in the energy norm[3] [35].

The Galerkin approximation $\boldsymbol{u}_h$ of $\boldsymbol{u}$ is obtained by replacing $\mathcal{V}^{\boldsymbol{g}}$ with the discrete space of admissible displacements given by

$$\mathcal{V}_h^{\boldsymbol{g}} = \mathcal{V}^{\boldsymbol{g}} \cap \mathcal{V}_h \tag{5}$$

in the minimization problem (1). Therefore,

$$\boldsymbol{u}_h = \underset{\boldsymbol{v} \in \mathcal{V}_h^{\boldsymbol{g}}}{\operatorname{argmin}} \left( \frac{1}{2} a(\boldsymbol{v}, \boldsymbol{v}) - f(\boldsymbol{v}) \right) \tag{6}$$

We are assuming here that the essential boundary conditions can be satisfied exactly in $\mathcal{V}_h$ since otherwise the intersection in (5) will be empty. In practice, the boundary data $\boldsymbol{g}$ is replaced by its nodal approximation $\boldsymbol{g}_h$ but the analysis of the effects of this error is classical [35] and so we shall ignore it to simplify the presentation.

Due to the conformity of $\mathcal{V}_h$, the strain energy associated with $\boldsymbol{v} \in \mathcal{V}_h$ is simply the sum of the contributions from the elements in the mesh. In other words,

$$a(\boldsymbol{v}, \boldsymbol{v}) = \sum_{E \in \mathcal{T}_h} a^E(\boldsymbol{v}, \boldsymbol{v}) \tag{7}$$

where we have denoted by $a^E$ the strain energy associated with element $E$ given by

$$a^E(\boldsymbol{u}, \boldsymbol{v}) = \int_E \boldsymbol{\sigma}(\boldsymbol{u}) : \boldsymbol{\epsilon}(\boldsymbol{v}) \mathrm{d}\boldsymbol{x} \tag{8}$$

As we will see, particular attention will be given in VEM to an appropriate approximation of these local strain energies that in turn determine the energetic behavior of the polyhedral elements. Naturally, this will first require a description of the element space $\mathcal{W}(E)$, which as discussed in the next section, will be a typical nodal finite element space.

[2]We refer to [16] for an example of set of mild geometric constraints that may be required of the elements in the mesh.

[3]The energy norm is equivalent to the $H^1$-norm in the space $\mathcal{V}^{\mathbf{0}}$.

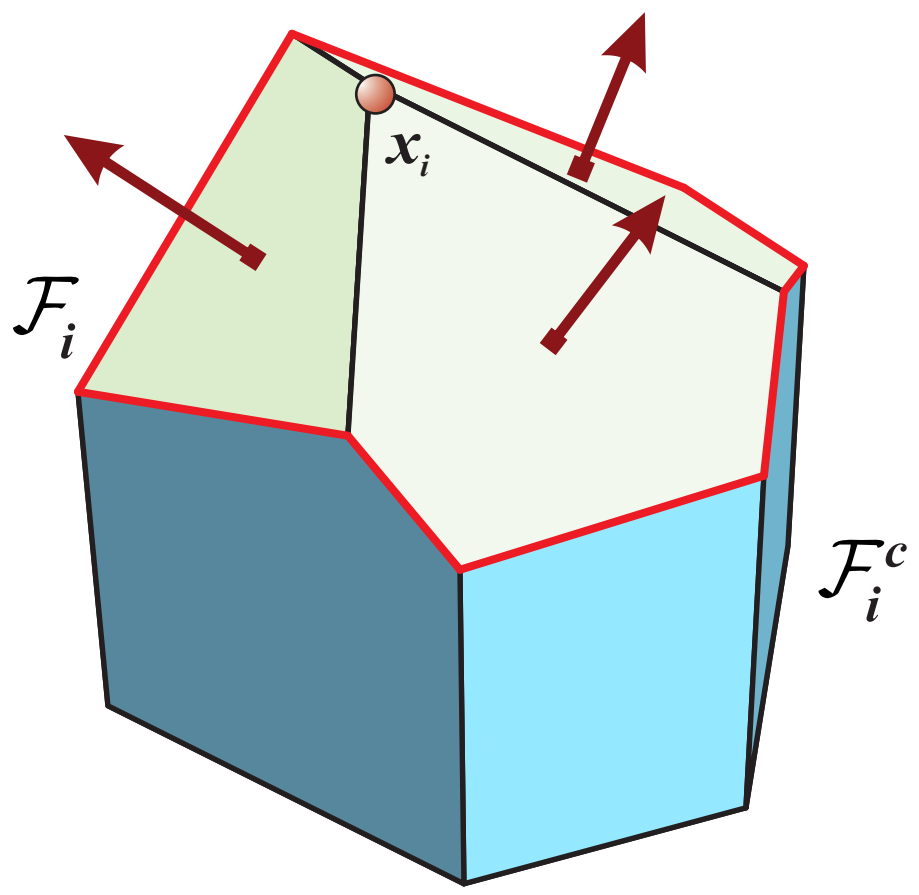

**Figure 2:** Decomposition of the faces of polyhedron into $\mathcal{F}_i$ and $\mathcal{F}_i^c$. Here, $\mathcal{F}_i$ represents the set of faces containing vertex $\boldsymbol{x}_i$ and $\mathcal{F}_i^c$ represents the remaining faces.

### *2.2. Construction and properties of* $\mathcal{W}(E)$

As stated before, each element $E \in \mathcal{T}_h$ is a polyhedron whose boundary consists of planar polygonal faces. Suppose $E$ has $n$ vertices located at $\boldsymbol{x}_1, \ldots, \boldsymbol{x}_n$. Let us denote by $\mathcal{F}_i$ the set of faces that include $\boldsymbol{x}_i$ and by $\mathcal{F}_i^c$ the remaining faces (see Fig. 2). Note that we do not require $E$ to be convex though convexity and its implications on the mesh topology can simplify certain aspects of implementation.

We will give a construction of the element space $\mathcal{W}(E)$ with three degrees of freedom associated with each vertex. To this effect, we consider the canonical basis $\boldsymbol{\varphi}_1, \ldots, \boldsymbol{\varphi}_{3n}$ of the form

$$\boldsymbol{\varphi}_{3i-2} = [\varphi_i, 0, 0]^\top, \quad \boldsymbol{\varphi}_{3i-1} = [0, \varphi_i, 0]^\top, \quad \boldsymbol{\varphi}_{3i} = [0, 0, \varphi_i]^\top, \qquad i = 1, \ldots, n \tag{9}$$

where $\varphi_1, \ldots, \varphi_n$ constitutes a set of *barycentric coordinates* for $E$. Examples of barycentric coordinates for polytopes can be found in [33, 12, 23, 25]. Among these, we will use the maximum entropy coordinates [23] later in our numerical studies. By definition, barycentric coordinates satisfy the Kronecker-delta property (i.e., $\varphi_i(\boldsymbol{x}_j) = \delta_{ij}$), which in turn implies that each $\boldsymbol{u} \in \mathcal{W}(E)$ is completely characterized by the values it assumes at the vertices of element $E$, consistent with the stated choice of degrees of freedom. Furthermore, $\varphi_i$ varies linearly along the edges of $E$ and vanishes on $\mathcal{F}_i^c$, i.e., the faces not incident on the associated vertex $\boldsymbol{x}_i$. *Also, the variation of* $\varphi_i$ *on the faces in* $\mathcal{F}_i$ *is determined* uniquely *by the geometry of those faces and independent of the shape of the element.* The latter two properties are crucial in guaranteeing inter-element continuity, and subsequently conformity of $\mathcal{V}_h$, as the variation of $\boldsymbol{u} \in \mathcal{W}(E)$ on a face is uniquely determined by values of $\boldsymbol{u}$ at vertices of that face and its geometry. Finally, barycentric coordinates can interpolate linear fields exactly, that is,

$$a + \boldsymbol{b} \cdot \boldsymbol{x} = \sum_{i=1}^{n} (a + \boldsymbol{b} \cdot \boldsymbol{x}_i)\, \varphi_i(\boldsymbol{x}) \tag{10}$$

for any $a \in \Re$ and $\boldsymbol{b} \in \Re^3$. This, in turn, implies that the element $E$ can represent rigid body motions and states of constant strain, i.e.,

$$\mathcal{W}(E) \supseteq \mathcal{P}(E) \doteq \left\{ \boldsymbol{a} + \boldsymbol{B}\boldsymbol{x} : \boldsymbol{a} \in \Re^3, \boldsymbol{B} \in \Re^{3\times 3} \right\} \tag{11}$$

guaranteeing the previously-stated first-order approximation capability of $\mathcal{W}(E)$. We should note that if $E$ is a tetrahedron, the well-known linear shape functions are the unique set of barycentric coordinates for $E$ and $\mathcal{W}(E) = \mathcal{P}(E)$.

As we will see in the next section, *what is more relevant in VEM is the behavior of functions in $\mathcal{W}(E)$ on the boundary, not in the interior of $E$.* Therefore, it is imperative to know the boundary behavior of the barycentric coordinate defining the basis functions (9). A useful observation is that, given any two-dimensional barycentric coordinates for planar polygons, we can use *harmonic lifting* to construct barycentric coordinates for $E$ exhibiting the above-mentioned properties. This process defines $\varphi_i$ as the solution to the Laplace equation whose boundary conditions are set to be the two-dimensional barycentric coordinates on faces in $\mathcal{F}_i$ and zero on the faces in $\mathcal{F}_i^c$ (see Fig. 2). The resulting coordinates $\varphi_i$'s will have the desired boundary behavior (Kronecker-delta property at the vertices, linearity and continuity along edges, and variation on the faces dictated by the choice of 2D barycentric coordinates), and their linear completeness follows from the linear completeness of 2D coordinates and properties of Laplace's equation. We remark that the use of *harmonic* basis, explicitly computed, have been explored in practice (see, for example, [14, 32]).

Later we shall assume the use of this harmonic construction of $\mathcal{W}(E)$ along with a particular set of boundary coordinates defined in [3], which possess the useful property that their average value can be computed explicitly based on the geometry of the underlying polygon (see the appendix for more details). If instead a nodal quadrature rule is used for computing surface integral encountered in the formulation, no distinction will be made between different barycentric coordinates (harmonic or otherwise) underlying $\mathcal{W}(E)$ (cf. Section 4).

We close this section by recalling the quasi-optimality of the error in the Galerkin solution, namely that the Galerkin error is bounded by a constant multiple of the error in the best approximation of $\boldsymbol{u}$ in $\mathcal{V}_h^{\boldsymbol{g}}$:

$$\|\boldsymbol{u} - \boldsymbol{u}_h\|_{1,\Omega} \leq C \inf_{\boldsymbol{w}_h \in \mathcal{V}_h^{\boldsymbol{g}}} \|\boldsymbol{u} - \boldsymbol{w}_h\|_{1,\Omega} \tag{12}$$

By the approximation property of the finite element space, the right hand side is $\mathcal{O}(h)$, and the convergence rate is linear.

## 3. Virtual Element Method (VEM)

While the Galerkin discretization on $\mathcal{T}_h$ is now completely defined, its realization is difficult to achieve in practice. The main source of difficulty is the evaluation of the weak form integrals, i.e., computing $a^E$ and $f$. Since the functions in $\mathcal{W}(E)$, and in particular its basis, are in general non-polynomial functions, available quadrature rules will inevitably lead to errors in the evaluation of the weak form integrals. Using high-order quadrature rules to reduce this error to acceptable levels is prohibitively expensive in practice since the construction of the basis functions (harmonics or otherwise), due to the lack of availability of explicit analytical expressions, is computational costly.

Acknowledging the presence of error in the evaluation of the linear and bilinear forms, we find ourselves committing a *variational crime* and deviating from the Galerkin framework, in effect replacing $a$, $a^E$ and $f$ by approximate mesh-dependent counterparts $a_h$, $a_h^E$ and $f_h$, respectively. The resulting approximate solution $\tilde{\boldsymbol{u}}_h$ minimizes the *discrete* total potential energy and is characterized by

$$\tilde{\boldsymbol{u}}_h = \underset{\boldsymbol{v} \in \mathcal{V}_h^{\boldsymbol{g}}}{\operatorname{argmin}} \left( \frac{1}{2} a_h(\boldsymbol{v}, \boldsymbol{v}) - f_h(\boldsymbol{v}) \right) \tag{13}$$

We can analyze the error $(\boldsymbol{u} - \tilde{\boldsymbol{u}}_h)$ by using the Galerkin solution $\boldsymbol{u}_h$ as an intermediary as follows:

$$\|\boldsymbol{u} - \tilde{\boldsymbol{u}}_h\|_{1,\Omega} \leq \|\boldsymbol{u} - \boldsymbol{u}_h\|_{1,\Omega} + \|\boldsymbol{u}_h - \tilde{\boldsymbol{u}}_h\|_{1,\Omega} \tag{14}$$

As before, the first term is governed by the approximation properties of the discrete space (cf. (12)). The second term represents the *consistency error* introduced by replacing the total potential energy by its discrete counterpart. According to Strang's lemma [20, 35], provided that $a_h$ is uniformly coercive on $\mathcal{V}_h^{\boldsymbol{0}}$, we have the following bound for this term:

$$\|\boldsymbol{u}_h - \tilde{\boldsymbol{u}}_h\|_{1,\Omega} \leq C \sup_{\boldsymbol{v} \in \mathcal{V}_h^{\boldsymbol{0}}} \frac{|a\left(\boldsymbol{u}_h, \boldsymbol{v}\right) - a_h\left(\boldsymbol{u}_h, \boldsymbol{v}\right)| + |f(\boldsymbol{v}) - f_h(\boldsymbol{v})|}{\|\boldsymbol{v}\|_{1,\Omega}} \tag{15}$$

Here, $C$ is a constant independent of $h$. If the discrete strain energy and load forms are defined such that the terms in the consistency error (15) are $\mathcal{O}(h)$, we can ensure that $\tilde{\boldsymbol{u}}_h$ converges to $\boldsymbol{u}$ at the same (optimal) rate as the Galerkin solution $\boldsymbol{u}_h$.

As before, within the conforming setting, the discrete energy bilinear form $a_h$ is typically obtained from the contribution of the discrete elemental bilinear forms

$$a_h(\boldsymbol{v}, \boldsymbol{v}) = \sum_{E \in \mathcal{T}_h} a_h^E(\boldsymbol{v}, \boldsymbol{v}) \tag{16}$$

As detailed in the remainder of this section, VEM gives a particular construction of $a_h^E$ such that $a_h^E(\boldsymbol{v}, \boldsymbol{v})$ and its variations are exact whenever $\boldsymbol{v}$ is either rigid body motion or a constant-strain displacement field on $E$. In other words, each element in the mesh will correctly represent the strain energy associated with these deformation states. This means that the so-called patch test will be passed at the element level. The consequence of the satisfaction of the element patch test, as discussed in Section 3.2, is that the consistency error introduced by replacing $a$ by $a_h$ (i.e., the first term in (15)) is $\mathcal{O}(h)$. Curiously, the construction of the discrete bilinear forms $a_E^h$ in VEM does not require numerical quadrature inside the element and therefore eliminates the need for costly computation of the basis functions in the interior of the element.

### *3.1. Kinematics decomposition of* $\mathcal{W}(E)$

We now discuss a particular kinematic decomposition of the deformation states in $\mathcal{W}(E)$ that is central to the VEM construction. In the remainder of this subsection, we will focus on element $E \in \mathcal{T}_h$ and thus omit the dependence on $E$ to ease the notation (e.g. write $\mathcal{W}$ for $\mathcal{W}(E)$ and $\mathcal{P}$ for $\mathcal{P}(E)$). For a function $\boldsymbol{w}$, we shall denote by $\overline{\boldsymbol{w}}$ the mean of the values it assumes over the vertices of $E$:

$$\overline{\boldsymbol{w}} = \frac{1}{n}\sum_{i=1}^{n} \boldsymbol{w}(\boldsymbol{x}_i) \tag{17}$$

This means, for example, that $\overline{\boldsymbol{x}}$ is the geometric center of $E$. Similarly, we will use $\langle \cdot \rangle$ to denote the volume average over $E$:

$$\langle \boldsymbol{w} \rangle = \frac{1}{|E|}\int_E \boldsymbol{w}\,\mathrm{d}\boldsymbol{x} \tag{18}$$

First let us split up $\mathcal{P}$, the space of linear displacements over $E$, into the spaces of rigid body motions and constant strain modes, defined, respectively, by

$$\mathcal{R} = \left\{\boldsymbol{a} + \boldsymbol{B}_{\mathsf{A}}\left(\boldsymbol{x} - \overline{\boldsymbol{x}}\right) : \boldsymbol{a} \in \Re^3, \boldsymbol{B}_{\mathsf{A}} \in \Re^{3\times 3}, \boldsymbol{B}_{\mathsf{A}}^{\top} = -\boldsymbol{B}_{\mathsf{A}}\right\} \tag{19}$$

$$\mathcal{C} = \left\{\boldsymbol{B}_{\mathsf{S}}\left(\boldsymbol{x} - \overline{\boldsymbol{x}}\right) : \boldsymbol{B}_{\mathsf{S}} \in \Re^{3\times 3}, \boldsymbol{B}_{\mathsf{S}}^{\top} = \boldsymbol{B}_{\mathsf{S}}\right\} \tag{20}$$

Observe that $\mathcal{P}$ is a direct sum of $\mathcal{R}$ and $\mathcal{C}$, and by (11), $\mathcal{R}$ and $\mathcal{C}$ are subspaces of $\mathcal{W}$. We next define bases for these spaces, respectively denoted by $\boldsymbol{r}_1, \ldots, \boldsymbol{r}_6$ and $\boldsymbol{c}_1, \ldots, \boldsymbol{c}_6$, as follows. We set $\boldsymbol{r}_1, \boldsymbol{r}_2, \boldsymbol{r}_3$ to be rigid body translation modes and $\boldsymbol{r}_4, \boldsymbol{r}_5, \boldsymbol{r}_6$ pure rotations about $\overline{\boldsymbol{x}}$:

$$\begin{aligned}
\boldsymbol{r}_1(\boldsymbol{x}) &= [1,0,0]^{\top} & \boldsymbol{r}_4(\boldsymbol{x}) &= \left[(\boldsymbol{x}-\overline{\boldsymbol{x}})_{(2)}, -(\boldsymbol{x}-\overline{\boldsymbol{x}})_{(1)}, 0\right]^{\top} \\
\boldsymbol{r}_2(\boldsymbol{x}) &= [0,1,0]^{\top} & \boldsymbol{r}_5(\boldsymbol{x}) &= \left[0, (\boldsymbol{x}-\overline{\boldsymbol{x}})_{(3)}, -(\boldsymbol{x}-\overline{\boldsymbol{x}})_{(2)}\right]^{\top} \\
\boldsymbol{r}_3(\boldsymbol{x}) &= [0,0,1]^{\top} & \boldsymbol{r}_6(\boldsymbol{x}) &= \left[-(\boldsymbol{x}-\overline{\boldsymbol{x}})_{(3)}, 0, (\boldsymbol{x}-\overline{\boldsymbol{x}})_{(1)}\right]^{\top}
\end{aligned} \tag{21}$$

Similarly, we choose $\boldsymbol{c}_1, \boldsymbol{c}_2, \boldsymbol{c}_3$ to correspond to deformations modes with constant axial strains, and $\boldsymbol{c}_4, \boldsymbol{c}_5, \boldsymbol{c}_6$ to represent three constant shear strains:

$$\begin{aligned}
\boldsymbol{c}_1(\boldsymbol{x}) &= \left[(\boldsymbol{x}-\overline{\boldsymbol{x}})_{(1)}, 0, 0\right]^{\top} & \boldsymbol{c}_4(\boldsymbol{x}) &= \left[(\boldsymbol{x}-\overline{\boldsymbol{x}})_{(2)}, (\boldsymbol{x}-\overline{\boldsymbol{x}})_{(1)}, 0\right]^{\top} \\
\boldsymbol{c}_2(\boldsymbol{x}) &= \left[0, (\boldsymbol{x}-\overline{\boldsymbol{x}})_{(2)}, 0\right]^{\top} & \boldsymbol{c}_5(\boldsymbol{x}) &= \left[0, (\boldsymbol{x}-\overline{\boldsymbol{x}})_{(3)}, (\boldsymbol{x}-\overline{\boldsymbol{x}})_{(2)}\right]^{\top} \\
\boldsymbol{c}_3(\boldsymbol{x}) &= \left[0, 0, (\boldsymbol{x}-\overline{\boldsymbol{x}})_{(3)}\right]^{\top} & \boldsymbol{c}_6(\boldsymbol{x}) &= \left[(\boldsymbol{x}-\overline{\boldsymbol{x}})_{(3)}, 0, (\boldsymbol{x}-\overline{\boldsymbol{x}})_{(1)}\right]^{\top}
\end{aligned} \tag{22}$$

In these expressions, the subscript within parentheses designates the component of the associated vector. The bases for $\mathcal{R}$ and $\mathcal{C}$ are illustrated for an arbitrary polyhedron in Fig. 3.

Next we define *projection* maps $\pi_{\mathcal{R}} : \mathcal{W} \to \mathcal{R}$ and $\pi_{\mathcal{C}} : \mathcal{W} \to \mathcal{C}$ that allow us to *extract* the rigid body motion and constant strain part of any deformation state $\boldsymbol{v} \in \mathcal{W}$. By definition, these maps will satisfy

$$\pi_{\mathcal{R}}\boldsymbol{r} = \boldsymbol{r}, \qquad \forall \boldsymbol{r} \in \mathcal{R} \tag{23}$$

$$\pi_{\mathcal{C}}\boldsymbol{c} = \boldsymbol{c}, \qquad \forall \boldsymbol{c} \in \mathcal{C} \tag{24}$$

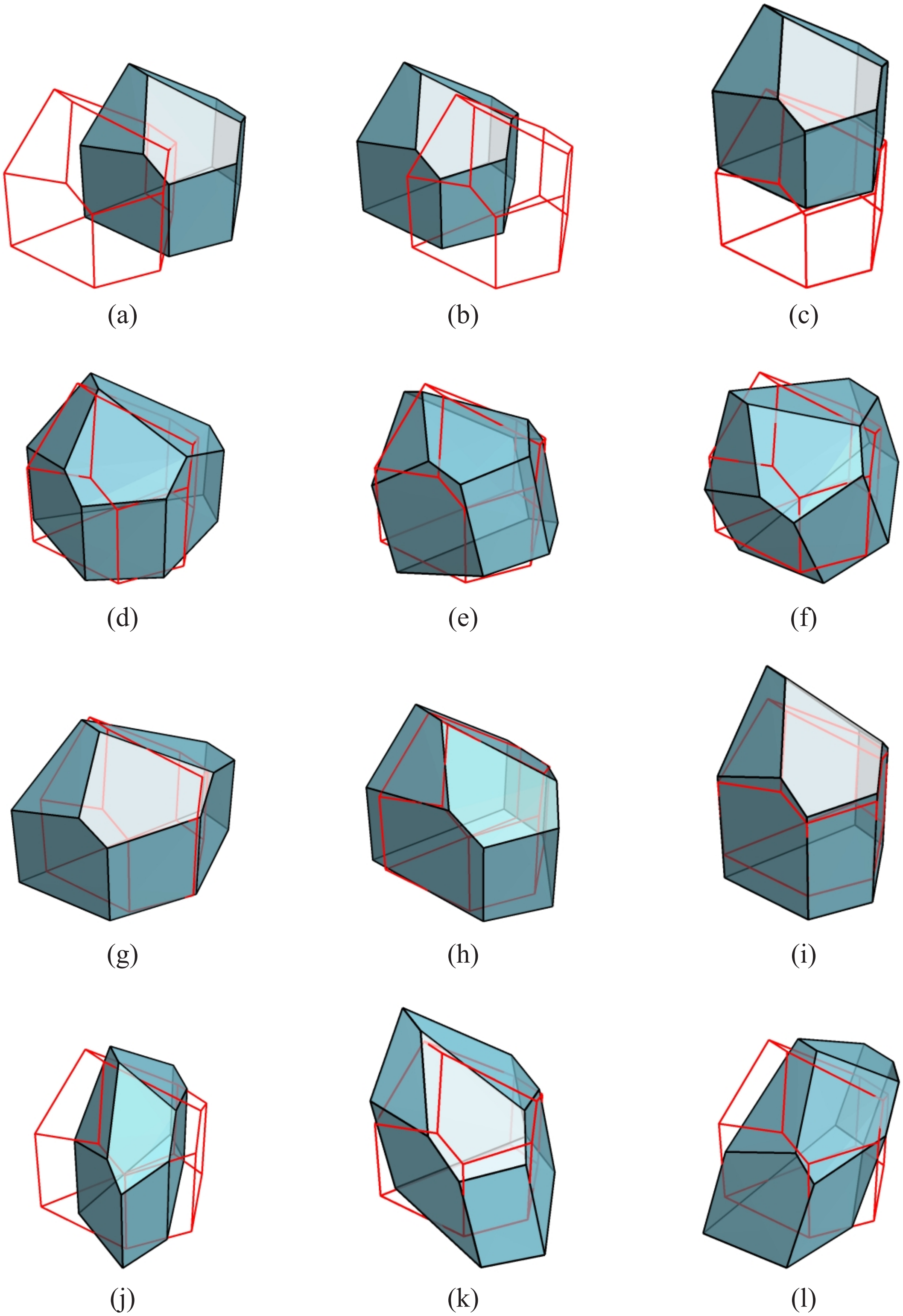

**Figure 3:** Illustration of the polynomial basis representing: rigid body translations (a) $\boldsymbol{r}_1$ (b) $\boldsymbol{r}_2$ (c) $\boldsymbol{r}_3$; rigid body rotations (d) $\boldsymbol{r}_4$ (e) $\boldsymbol{r}_5$ (f) $\boldsymbol{r}_6$; axial strains (g) $\boldsymbol{c}_1$ (h) $\boldsymbol{c}_2$ (i) $\boldsymbol{c}_3$; shear strains (j) $\boldsymbol{c}_4$ (k) $\boldsymbol{c}_5$ (l) $\boldsymbol{c}_6$.

Moreover, it will be useful to require the following *orthogonality* conditions:

$$\pi_{\mathcal{R}}\boldsymbol{c} = \boldsymbol{0}, \qquad \forall \boldsymbol{c} \in \mathcal{C} \tag{25}$$

$$\pi_{\mathcal{C}}\boldsymbol{r} = \boldsymbol{0}, \qquad \forall \boldsymbol{r} \in \mathcal{R} \tag{26}$$

reflecting the fact that elements of $\mathcal{C}$ will contain no rigid body motions and, similarly, elements of $\mathcal{R}$ will be associated with null element in $\mathcal{C}$. Subsequently, $\pi_{\mathcal{C}}\pi_{\mathcal{R}} = \pi_{\mathcal{R}}\pi_{\mathcal{C}} = 0$ and

$$\pi_{\mathcal{P}} = \pi_{\mathcal{R}} + \pi_{\mathcal{C}} \tag{27}$$

defines a projection onto $\mathcal{P}$.

A projection map $\pi_{\mathcal{R}}$ satisfying the above properties is given by

$$\pi_{\mathcal{R}}\boldsymbol{v} = \overline{\boldsymbol{v}} + \langle \boldsymbol{\omega}(\boldsymbol{v}) \rangle (\boldsymbol{x} - \overline{\boldsymbol{x}}) \tag{28}$$

Recall that $\boldsymbol{\omega}(\cdot)$ is the skew-symmetric gradient operator. Observe that we defined the space $\mathcal{C}$ such that $\overline{\boldsymbol{c}} = \boldsymbol{0}$ and $\boldsymbol{\omega}(\boldsymbol{c}) = \boldsymbol{0}$ for all $\boldsymbol{c} \in \mathcal{C}$ and so (25) immediately follows from the definition of $\pi_{\mathcal{R}}$. It is also straightforward to verify (23) since for $\boldsymbol{v} = \boldsymbol{a} + \boldsymbol{B}_{\mathsf{A}}(\boldsymbol{x} - \overline{\boldsymbol{x}})$, with $\boldsymbol{B}_{\mathsf{A}}$ an antisymmetric tensor: we have $\overline{\boldsymbol{v}} = \boldsymbol{a}$ and $\langle \boldsymbol{\omega}(\boldsymbol{v}) \rangle = \left\langle \boldsymbol{B}_{\mathsf{A}} - \boldsymbol{B}_{\mathsf{A}}^{\top} \right\rangle / 2 = \boldsymbol{B}_{\mathsf{A}}$ and so, by definition, $\pi_{\mathcal{R}}\boldsymbol{v} = \boldsymbol{v}$.

We also note that

$$\overline{\pi_{\mathcal{R}}\boldsymbol{v}} = \overline{\boldsymbol{v}}, \qquad \boldsymbol{\omega}(\pi_{\mathcal{R}}\boldsymbol{v}) = \langle \boldsymbol{\omega}(\boldsymbol{v}) \rangle, \qquad \boldsymbol{\epsilon}(\pi_{\mathcal{R}}\boldsymbol{v}) = \boldsymbol{0} \tag{29}$$

The first two relations show that the translation and rotation of $\pi_{\mathcal{R}}\boldsymbol{v}$ are equal to *average* translation and rotation of $\boldsymbol{v}$, respectively. Finally, it will be useful (later) to express $\pi_{\mathcal{R}}$ in terms of the basis of $\mathcal{R}$ as follows

$$\pi_{\mathcal{R}}\boldsymbol{v} = (\overline{\boldsymbol{v}})_{(1)}\boldsymbol{r}_1 + (\overline{\boldsymbol{v}})_{(2)}\boldsymbol{r}_2 + (\overline{\boldsymbol{v}})_{(3)}\boldsymbol{r}_3 + \langle \boldsymbol{\omega}(\boldsymbol{v}) \rangle_{(12)}\boldsymbol{r}_4 + \langle \boldsymbol{\omega}(\boldsymbol{v}) \rangle_{(23)}\boldsymbol{r}_5 + \langle \boldsymbol{\omega}(\boldsymbol{v}) \rangle_{(31)}\boldsymbol{r}_6 \tag{30}$$

One important observation is the volumetric integral in the definition of $\pi_{\mathcal{R}}$ can be transformed as a boundary integral since

$$\int_E \boldsymbol{\omega}(\boldsymbol{v})\mathrm{d}\boldsymbol{x} = \frac{1}{2}\int_E \left(\nabla\boldsymbol{v} - \nabla^{\top}\boldsymbol{v}\right)\mathrm{d}\boldsymbol{x} = \frac{1}{2}\int_{\partial E} (\boldsymbol{v} \otimes \boldsymbol{n} - \boldsymbol{n} \otimes \boldsymbol{v})\,\mathrm{d}\boldsymbol{s} \tag{31}$$

where $\boldsymbol{n}$ is the unit normal to $\partial E$ pointing outwards. Hence, $\pi_{\mathcal{R}}\boldsymbol{v}$ is completely determined by the nodal (vertex) values of $\boldsymbol{v}$. Recall that over each face of $E$, the variation of $\boldsymbol{v} \in \mathcal{W}$ is completely determined by its values at the vertices of that face.

We define the projection onto the space of constant-strain fields as

$$\pi_{\mathcal{C}}\boldsymbol{v} = \langle \boldsymbol{\epsilon}(\boldsymbol{v}) \rangle (\boldsymbol{x} - \overline{\boldsymbol{x}}) \tag{32}$$

One can again directly verify conditions (24) and (26) from this definition. Analogous to (29), we have

$$\overline{\pi_{\mathcal{C}}\boldsymbol{v}} = \boldsymbol{0}, \qquad \boldsymbol{\omega}(\pi_{\mathcal{C}}\boldsymbol{v}) = \boldsymbol{0}, \qquad \boldsymbol{\epsilon}(\pi_{\mathcal{C}}\boldsymbol{v}) = \langle \boldsymbol{\epsilon}(\boldsymbol{v}) \rangle \tag{33}$$

illustrating the strain associated with $\pi_{\mathcal{C}}\boldsymbol{v}$ is the average of strain of $\boldsymbol{v}$, and $\pi_{\mathcal{C}}\boldsymbol{v}$ does not contain any rigid body translation or rotation. We can also express the expansion of $\pi_{\mathcal{C}}$ in terms of the basis for $\mathcal{C}$ as follows:

$$\pi_{\mathcal{C}}\boldsymbol{v} = \langle\boldsymbol{\epsilon}(\boldsymbol{v})\rangle_{(11)}\,\boldsymbol{c}_1 + \langle\boldsymbol{\epsilon}(\boldsymbol{v})\rangle_{(22)}\,\boldsymbol{c}_2 + \langle\boldsymbol{\epsilon}(\boldsymbol{v})\rangle_{(33)}\,\boldsymbol{c}_3 + \langle\boldsymbol{\epsilon}(\boldsymbol{v})\rangle_{(12)}\,\boldsymbol{c}_4 + \langle\boldsymbol{\epsilon}(\boldsymbol{v})\rangle_{(23)}\,\boldsymbol{c}_5 + \langle\boldsymbol{\epsilon}(\boldsymbol{v})\rangle_{(31)}\,\boldsymbol{c}_6 \tag{34}$$

Finally, as with $\pi_{\mathcal{R}}$, the volume integral in $\pi_{\mathcal{C}}\boldsymbol{v}$ can be written as an integral over the boundary of $E$

$$\int_E \boldsymbol{\epsilon}(\boldsymbol{v})\mathrm{d}\boldsymbol{x} = \frac{1}{2}\int_E \left(\nabla\boldsymbol{v} + \nabla^\top\boldsymbol{v}\right)\mathrm{d}\boldsymbol{x} = \frac{1}{2}\int_{\partial E}\left(\boldsymbol{v}\otimes\boldsymbol{n} + \boldsymbol{n}\otimes\boldsymbol{v}\right)\mathrm{d}\boldsymbol{s} \tag{35}$$

and so $\pi_{\mathcal{C}}\boldsymbol{v}$ is determined by the nodal values of $\boldsymbol{v}$.

An important property of $\pi_{\mathcal{C}}$ is that for all $\boldsymbol{v}\in\mathcal{W}$, the term $\boldsymbol{v}-\pi_{\mathcal{C}}\boldsymbol{v}$ is *energetically* orthogonal to $\mathcal{C}$, that is

$$a^E(\boldsymbol{c}, \boldsymbol{v}-\pi_{\mathcal{C}}\boldsymbol{v}) = 0, \quad \forall \boldsymbol{c}\in\mathcal{C} \tag{36}$$

As we will see in the next section, *this property plays a crucial role in ensuring consistency of the VEM bilinear form.* To verify this identity, we appeal to the last equality in (33) and the fact that $\boldsymbol{\sigma}(\boldsymbol{c})$ is a constant field:

$$a^E(\boldsymbol{c}, \boldsymbol{v}-\pi_{\mathcal{C}}\boldsymbol{v}) = \int_E \boldsymbol{\sigma}(\boldsymbol{c}) : \left[\boldsymbol{\epsilon}(\boldsymbol{v}) - \boldsymbol{\epsilon}(\pi_{\mathcal{C}}\boldsymbol{v})\right]\mathrm{d}\boldsymbol{x} = \boldsymbol{\sigma}(\boldsymbol{c}) : \left[\int_E \boldsymbol{\epsilon}(\boldsymbol{v})\mathrm{d}\boldsymbol{x} - \boldsymbol{\epsilon}(\pi_{\mathcal{C}}\boldsymbol{v})\,|E|\right] = 0 \tag{37}$$

In fact, we can show that (32) is the only projection onto that satisfies (36). *In other words, the energy orthogonality condition uniquely determines the projection on* $\mathcal{C}$[4].

As rigid body motions have zero strain and thus contain no strain energy, we can see that the energy orthogonality extends to $\mathcal{P}$ and $\pi_{\mathcal{P}}$ (cf. equation (27)) as

$$a^E(\boldsymbol{p}, \boldsymbol{v}-\pi_{\mathcal{P}}\boldsymbol{v}) = 0, \quad \forall \boldsymbol{p}\in\mathcal{P} \tag{38}$$

for all $\boldsymbol{v}\in\mathcal{W}$.

With the projection maps defined explicitly, we can obtain an additive kinematic decomposition of a given deformation state $\boldsymbol{v}\in\mathcal{W}$ into its rigid body motion, constant strain and the remaining *higher-order* components:

$$\boldsymbol{v} = \pi_{\mathcal{R}}\boldsymbol{v} + \pi_{\mathcal{C}}\boldsymbol{v} + (\boldsymbol{v}-\pi_{\mathcal{P}}\boldsymbol{v}) \tag{39}$$

The remainder $\boldsymbol{v}-\pi_{\mathcal{P}}\boldsymbol{v}$ belongs to a $(3n-12)$-dimensional subspace of $\mathcal{W}$, which we shall denote by $\mathcal{H}$. This space consists of displacements modes that are either higher-order polynomials or non-polynomial functions. For example, if $E$ is a cube, and $\mathcal{W}$ is the space of trilinear displacement fields, $\boldsymbol{v}-\pi_{\mathcal{P}}\boldsymbol{v}$ is a linear combination of 12 *hourglass* modes consisting of high-order polynomials. For a distorted hexahedron, these modes will consists of rational functions even with the classical iso-parametric finite element bases [13].

[4] Observe, for instance, that (36) immediately implies (26) since $a^E(\pi_{\mathcal{C}}\boldsymbol{r}, \pi_{\mathcal{C}}\boldsymbol{r}) = a^E(\boldsymbol{r}, \pi_{\mathcal{C}}\boldsymbol{r}) = 0$ for all $\boldsymbol{r}\in\mathcal{R}$.

### 3.2. Construction of the discrete bilinear forms

The kinematic decomposition of deformation state $\boldsymbol{v}$ in the form (39), by the virtue of energy-orthogonality condition (36), leads to a decomposition of strain energy associated with $\boldsymbol{v}$:

$$\begin{aligned} a^E(\boldsymbol{v},\boldsymbol{v}) &= a^E\left(\pi_{\mathcal{R}}\boldsymbol{v}+\pi_{\mathcal{C}}\boldsymbol{v}+(\boldsymbol{v}-\pi_{\mathcal{P}}\boldsymbol{v}),\pi_{\mathcal{R}}\boldsymbol{v}+\pi_{\mathcal{C}}\boldsymbol{v}+(\boldsymbol{v}-\pi_{\mathcal{P}}\boldsymbol{v})\right) \\ &= a^E(\pi_{\mathcal{C}}\boldsymbol{v},\pi_{\mathcal{C}}\boldsymbol{v})+2a^E(\pi_{\mathcal{C}}\boldsymbol{v},\boldsymbol{v}-\pi_{\mathcal{P}}\boldsymbol{v})+a^E(\boldsymbol{v}-\pi_{\mathcal{P}}\boldsymbol{v},\boldsymbol{v}-\pi_{\mathcal{P}}\boldsymbol{v}) \\ &= a^E(\pi_{\mathcal{C}}\boldsymbol{v},\pi_{\mathcal{C}}\boldsymbol{v})+a^E(\boldsymbol{v}-\pi_{\mathcal{P}}\boldsymbol{v},\boldsymbol{v}-\pi_{\mathcal{P}}\boldsymbol{v}) \end{aligned} \tag{40}$$

In the first equality, we have used the linearity and symmetry of $a^E$ along with the fact that $\boldsymbol{\epsilon}(\pi_{\mathcal{R}}\boldsymbol{v})=\boldsymbol{0}$. The identity (38) and $\pi_{\mathcal{C}}\boldsymbol{v}\in\mathcal{P}$, are used in the second equality to eliminate the coupling term. The first term in the right hand side of the final expression in (40) represents the energy associated with the constant-strain component of $\boldsymbol{v}$ while the second term gives the energy associated with the remaining higher-order part. Observe that the first term can be computed exactly with the knowledge of the volume of $E$ since its integrand, $\boldsymbol{\sigma}(\pi_{\mathcal{C}}\boldsymbol{v}):\boldsymbol{\epsilon}(\pi_{\mathcal{C}}\boldsymbol{v})$, is a constant field.

Now comes another key observation: we can replace the second term in (40) by a crude estimate, one that can be conveniently computed, without affecting the energy associated with rigid body motion and constant-strain component of $\boldsymbol{v}$. This suggests defining the following discrete energy form for $E$:

$$a_h^E(\boldsymbol{u},\boldsymbol{v})\doteq a^E(\pi_{\mathcal{C}}\boldsymbol{u},\pi_{\mathcal{C}}\boldsymbol{v})+s^E(\boldsymbol{u}-\pi_{\mathcal{P}}\boldsymbol{u},\boldsymbol{v}-\pi_{\mathcal{P}}\boldsymbol{v}) \tag{41}$$

where $s^E$ is a prescribed symmetric continuous bilinear form on $\mathcal{W}$.

Noting that for $\boldsymbol{h}\in\mathcal{H}$,

$$a_h^E(\boldsymbol{h},\boldsymbol{h})=s^E(\boldsymbol{h},\boldsymbol{h}) \tag{42}$$

it is evident that $s^E$ must be positive definite on the space of higher-order deformations $\mathcal{H}$. Otherwise, non-zero higher-order deformation modes may be assigned zero strain energy, potentially leading to global zero-energy modes and rank deficiency of the global system. In general, we may not have a guarantee of uniform coercivity of $a_h$ which is required for establishing estimate (15). A suitable choice of $s^E$ thus ensures the *stability* of the method by guaranteeing that discrete bilinear form inherits the coercivity of the exact bilinear form[5]. As mentioned before, for hexahedral elements, $\mathcal{H}$ is the space of hourglass modes, and so, in light of (42), $s^E$ essentially prescribes the "hourglass" stiffness (see, for example, [13]).

By virtue of the decomposition, the choice of $s^E$ does not affect the *polynomial consistency* of $a_h^E$. Indeed, for $\boldsymbol{p}\in\mathcal{P}$ and $\boldsymbol{v}\in\mathcal{W}$, we have

$$\begin{aligned} a_h^E(\boldsymbol{p},\boldsymbol{v}) &= a^E(\pi_{\mathcal{C}}\boldsymbol{p},\pi_{\mathcal{C}}\boldsymbol{v})+s^E(\boldsymbol{p}-\pi_{\mathcal{P}}\boldsymbol{p},\boldsymbol{v}-\pi_{\mathcal{P}}\boldsymbol{v}) \\ &= a^E(\pi_{\mathcal{C}}\boldsymbol{p},\pi_{\mathcal{C}}\boldsymbol{v}) && \text{(since } \boldsymbol{p}-\pi_{\mathcal{P}}\boldsymbol{p}=\boldsymbol{0}) \\ &= a^E(\pi_{\mathcal{C}}\boldsymbol{p},\boldsymbol{v}) && \text{(by (36))} \\ &= a^E(\boldsymbol{p},\boldsymbol{v}) && \text{(since } \boldsymbol{\epsilon}(\pi_{\mathcal{C}}\boldsymbol{p})=\boldsymbol{\epsilon}(\pi_{\mathcal{P}}\boldsymbol{p})=\boldsymbol{\epsilon}(\boldsymbol{p})) \end{aligned} \tag{43}$$

[5]A sufficient condition for this is that for some positive constants $\beta_1$ and $\beta_2$, independent of $h$ and $E$, we have $\beta_1 a^E(\boldsymbol{h},\boldsymbol{h})\le s^E(\boldsymbol{h},\boldsymbol{h})\le\beta_2 a^E(\boldsymbol{h},\boldsymbol{h})$ for all $\boldsymbol{h}\in\mathcal{H}$. This means that the strain energy associated with higher-order modes, as prescribed by $s^E$, scale uniformly with the exact strain energy.

This shows that the discrete bilinear form will exactly capture the strain energy associated with the linear deformation $\boldsymbol{p}$ and its variations. As a result, an element based on (41) will pass the first order patch test. The main consequence of this is that the error introduced by replacing the exact strain energy with $a_h$, the first term in (15), is $\mathcal{O}(h)$. To see this, we first split up exact and discrete strain energies into the element contributions (cf. (7) and (16)) and use (43) to add and subtract a term containing $\pi_{\mathcal{P}(E)}\boldsymbol{u}_h$ over each element:

$$\begin{aligned}
\left|a\left(\boldsymbol{u}_h,\boldsymbol{v}\right)-a_h\left(\boldsymbol{u}_h,\boldsymbol{v}\right)\right| &\leq \sum_{E\in\mathcal{T}_h}\left|a^E\left(\boldsymbol{u}_h,\boldsymbol{v}\right)-a_h^E\left(\boldsymbol{u}_h,\boldsymbol{v}\right)\right| \\
&= \sum_{E\in\mathcal{T}_h}\left|a^E\left(\boldsymbol{u}_h,\boldsymbol{v}\right)-a^E(\pi_{\mathcal{P}(E)}\boldsymbol{u}_h,\boldsymbol{v})+a_h^E(\pi_{\mathcal{P}(E)}\boldsymbol{u}_h,\boldsymbol{v})-a_h^E\left(\boldsymbol{u}_h,\boldsymbol{v}\right)\right| \\
&= \sum_{E\in\mathcal{T}_h}\left|a^E\left(\boldsymbol{u}_h-\pi_{\mathcal{P}(E)}\boldsymbol{u}_h,\boldsymbol{v}\right)-a_h^E\left(\boldsymbol{u}_h-\pi_{\mathcal{P}(E)}\boldsymbol{u}_h,\boldsymbol{v}\right)\right| \qquad (44) \\
&\leq C\sum_{E\in\mathcal{T}_h}\left\|\boldsymbol{u}_h-\pi_{\mathcal{P}(E)}\boldsymbol{u}_h\right\|_{1,E}\|\boldsymbol{v}\|_{1,E} \qquad \text{(using continuity of } a^E, a_h^E\text{)} \\
&\leq C\left(\sum_{E\in\mathcal{T}_h}\left\|\boldsymbol{u}_h-\pi_{\mathcal{P}(E)}\boldsymbol{u}_h\right\|_{1,E}^2\right)^{1/2}\left(\sum_{E\in\mathcal{T}_h}\|\boldsymbol{v}\|_{1,E}^2\right)^{1/2} \\
&\leq C'h\|\boldsymbol{v}\|_{1,\Omega} \qquad (45)
\end{aligned}$$

Here $C$ and $C'$ are constants independent of the mesh size. The last inequality, roughly speaking, is a consequence of the fact that $\pi_{\mathcal{P}(E)}\boldsymbol{u}_h$ is a first-order approximation to $\boldsymbol{u}_h$ and assumes uniform shape-regularity of the mesh.

Regarding choice of the bilinear form $s^E$, there is quite a bit of freedom in practice since $s^E$ can be any approximation of exact strain energy $a^E$, so long as it respects the stability requirement. On one end of the spectrum, we can define $s^E$ through quadrature as

$$s^E(\boldsymbol{u},\boldsymbol{v}) \doteq ⨎_E \boldsymbol{\sigma}(\boldsymbol{u}):\boldsymbol{\epsilon}(\boldsymbol{v})\mathrm{d}\boldsymbol{x} \qquad (46)$$

Here $⨎_E$ indicates that the volume integral is evaluated using a suitable quadrature rule. This approach has been pursued in [37] in order to alleviate the burden of numerical quadrature in the finite element setting. A simple choice on the other end is given by [10]:

$$s^E(\boldsymbol{u},\boldsymbol{v}) = \sum_{i=1}^{n}\alpha^E\boldsymbol{u}(\boldsymbol{x}_i)\cdot\boldsymbol{v}(\boldsymbol{x}_i) \qquad (47)$$

where $\alpha^E$ is a positive parameter that ensures the right scale of strain energies assigned to higher-order modes. Determining an appropriate value of this constant is deferred to the next section. Observe that this definition only involves the nodal values of the arguments and therefore eliminates the need for volumetric quadrature and construction of the basis functions inside the element. Such a choice highlights the advantage of VEM in constructing an inexpensive discretization scheme on arbitrary meshes.

We remark that the flexibility in the choice of $s^E$ can be exploited to enhance other characteristics of the resulting method (e.g. satisfaction of a discrete maximum principle or improvement of performance of algebraic solvers) as illustrated in [27, 28].

### *3.3. Construction of* $f_h(\boldsymbol{v})$

It is also possible to construct a first-order accurate discrete load linear form $f_h$ without the need for basis functions in the interior of the elements. We can do so by defining a *nodal* quadrature scheme over each element to treat the body force term and a *nodal* quadrature scheme over each face in the traction boundary $\Gamma_{\boldsymbol{t}}$. For each $\boldsymbol{v} \in \mathcal{V}_h$, we thus set

$$f_h(\boldsymbol{v}) \doteq \sum_{E \in \mathcal{T}_h} \fint_E \boldsymbol{b} \cdot \boldsymbol{v} \mathrm{d}\boldsymbol{x} + \sum_{F \subseteq \Gamma_{\boldsymbol{t}}} \fint_F \boldsymbol{t} \cdot \boldsymbol{v} \mathrm{d}\boldsymbol{s} \tag{48}$$

with two quadrature schemes[6] defined as follows (cf. Fig. 4a).

The two-dimensional surface integral over a face $F$ with $m$ vertices is approximated by

$$\fint_F T \mathrm{d}\boldsymbol{s} \doteq \sum_{j=1}^{m} w_j^F T(\boldsymbol{x}_j^F) \tag{49}$$

Here $\boldsymbol{x}_1^F, \ldots, \boldsymbol{x}_m^F$ denote the location of vertices incident on face $F$. The weight $w_j^F$ is the area of quadrilateral formed by $\boldsymbol{x}_j^F$, the midpoint of edges incident on $\boldsymbol{x}_j^F$, and the centroid of $F$. It is clear that $\sum_{j=1}^m w_j^F = |F|$ so the quadrature can integrate constant functions exactly. It is straightforward to show that this scheme can also exactly integrate linear fields.

The nodal quadrature for volume integral over $E$ is defined in a similar manner as:

$$\fint_E K \mathrm{d}\boldsymbol{x} \doteq \sum_{i=1}^{n} w_i^E K(\boldsymbol{x}_i) \tag{50}$$

where $\boldsymbol{x}_1, \ldots, \boldsymbol{x}_n$ denote the location of vertices of $E$, and the weight $w_i^E$ is the volume of polyhedron formed by $\boldsymbol{x}_i$, the centroid of $E$, the centroid of faces in $\mathcal{F}_i$, and the centroid of edges incident on $\boldsymbol{x}_i$. Again, we can show that $\fint_E p \mathrm{d}\boldsymbol{x} = \int_E p \mathrm{d}\boldsymbol{x}$, when $p$ is a linear polynomial.

The consequence of first-order accuracy of these quadrature schemes is that, when $\boldsymbol{b}$ and $\boldsymbol{t}$ are sufficiently smooth fields, the error in replacing the exact linear form $f$ with the approximation $f_h$, i.e., the second term in 15, is $\mathcal{O}(h)$. Therefore, convergence rate remains to be linear. Finally, we observe that the both quadrature schemes, by virtue of being nodal, do not distinguish between the choice of barycentric coordinates used in construction of element spaces $\mathcal{W}(E)$.

## 4. Implementation aspects

In this section, we discuss in detail the implementation of the VEM formulation presented in the previous section. More specifically, we will give explicit expressions for the element stiffness matrix corresponding to the discrete bilinear form defined in (41). The global stiffness matrix is obtained via the standard finite element assembly process (corresponding to the summation in (16)). Our derivations involve first obtaining discrete representations of the projection maps. At the end of this section, comments will be made on the structure of the element stiffness matrices which will further elucidates the nature of the approximation of element strain energies in VEM. Throughout this section, we will again focus our attention to an element $E$ in the mesh with vertices located at $\boldsymbol{x}_1, \ldots, \boldsymbol{x}_n$.

[6] Such quadrature rules have been used in the MFD literature see for example, [19, 9].

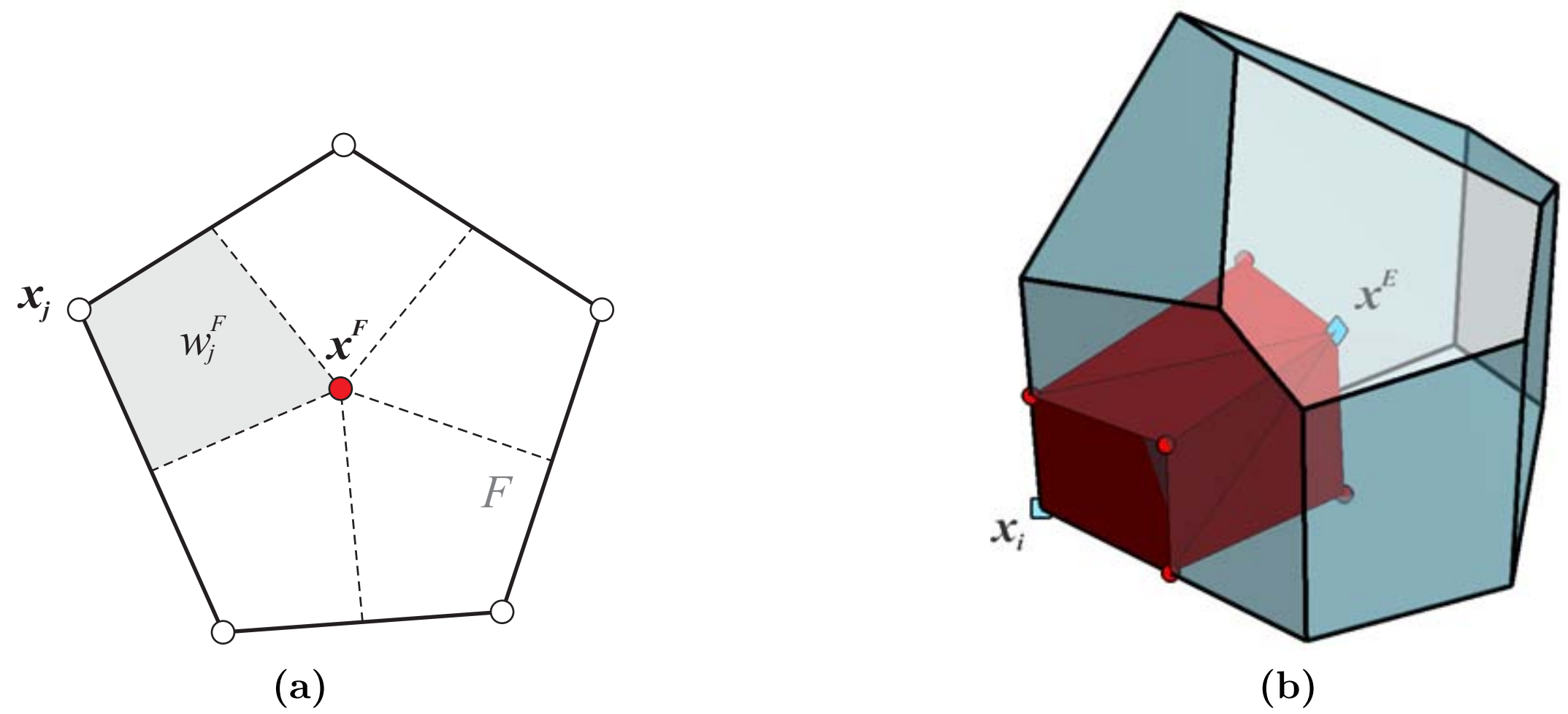

**Figure 4:** Illustration of nodal quadrature schemes for (a) surface integration: $\boldsymbol{x}^F$ and $w_j^F$ represent the centroid and the nodal weight associated with the vertex $\boldsymbol{x}_j$ on the face $F$, respectively; (b) volume integration: the nodal weight corresponding to the vertex $\boldsymbol{x}_i$ is volume of the red polyhedron formed by element centroid $\boldsymbol{x}^E$, three face centroids, three edge mid-points and the vertex $\boldsymbol{x}_i$. It is assumed that $F$ and $E$ are *star-shaped* with respect to $\boldsymbol{x}^F$ and $\boldsymbol{x}^E$, respectively.

### *4.1. Expressions for projection matrices*

Let us first define the map $\boldsymbol{\chi} : \mathcal{W} \to \mathfrak{R}^{3n}$ that computes nodal values (i.e., the degrees of freedom) associated with deformation states over $E$. For $i = 1, \ldots, n$,

$$[\boldsymbol{\chi}(\boldsymbol{v})]_{(3i-2)} = [\boldsymbol{v}(\boldsymbol{x}_i)]_{(1)}\,, \quad [\boldsymbol{\chi}(\boldsymbol{v})]_{(3i-1)} = [\boldsymbol{v}(\boldsymbol{x}_i)]_{(2)}\,, \quad [\boldsymbol{\chi}(\boldsymbol{v})]_{(3i)} = [\boldsymbol{v}(\boldsymbol{x}_i)]_{(3)} \tag{51}$$

and thus we have the expansion:

$$\boldsymbol{v}(\boldsymbol{x}) = \sum_{j=1}^{3n} [\boldsymbol{\chi}(\boldsymbol{v})]_{(j)}\, \boldsymbol{\varphi}_j(\boldsymbol{x}) \tag{52}$$

Note that any $\boldsymbol{v} \in \mathcal{W}$ can be uniquely identified with $\boldsymbol{\chi}(\boldsymbol{v}) \in \mathfrak{R}^{3n}$, and conversely, any array in $\mathfrak{R}^{3n}$ identifies with a member of $\mathcal{W}$.

The discrete counterpart to projection maps $\pi_{\mathcal{R}}, \pi_{\mathcal{C}}$ are $3n \times 3n$ matrices $\boldsymbol{P}_{\mathcal{R}}, \boldsymbol{P}_{\mathcal{C}}$ whose application to the nodal representation of $\boldsymbol{v}$ gives the nodal representation of its projection, that is,

$$\boldsymbol{P}_{\mathcal{R}}\boldsymbol{\chi}(\boldsymbol{v}) = \boldsymbol{\chi}(\pi_{\mathcal{R}}\boldsymbol{v})\,, \qquad \boldsymbol{P}_{\mathcal{C}}\boldsymbol{\chi}(\boldsymbol{v}) = \boldsymbol{\chi}(\pi_{\mathcal{C}}\boldsymbol{v}) \tag{53}$$

Setting $\boldsymbol{v} = \boldsymbol{\varphi}_j$ in the above expressions yields the alternative characterization of these matrices as

$$\pi_{\mathcal{R}}\boldsymbol{\varphi}_j = \sum_{k=1}^{3n} (\boldsymbol{P}_{\mathcal{R}})_{(kj)}\, \boldsymbol{\varphi}_k, \qquad \pi_{\mathcal{C}}\boldsymbol{\varphi}_j = \sum_{k=1}^{3n} (\boldsymbol{P}_{\mathcal{C}})_{(kj)}\, \boldsymbol{\varphi}_k \tag{54}$$

To obtain an explicit expression for $\boldsymbol{P}_{\mathcal{R}}$, we use (30) to write,

$$\pi_{\mathcal{R}}\boldsymbol{\varphi}_j = \sum_{\ell=1}^{6} (\boldsymbol{W}_{\mathcal{R}})_{(j\ell)}\, \boldsymbol{r}_\ell \tag{55}$$

Here, $\boldsymbol{W}_{\mathcal{R}}$ is the $3n \times 6$ matrix whose $j$th row is given by

$$\left[\left(\overline{\boldsymbol{\varphi}_j}\right)_{(1)}, \left(\overline{\boldsymbol{\varphi}_j}\right)_{(2)}, \left(\overline{\boldsymbol{\varphi}_j}\right)_{(3)}, \left\langle \boldsymbol{\omega}(\boldsymbol{\varphi}_j)\right\rangle_{(12)}, \left\langle \boldsymbol{\omega}(\boldsymbol{\varphi}_j)\right\rangle_{(23)}, \left\langle \boldsymbol{\omega}(\boldsymbol{\varphi}_j)\right\rangle_{(31)}\right] \tag{56}$$

Expanding $\boldsymbol{r}_\ell$ in terms of the canonical basis functions, $\boldsymbol{\varphi}_k$, we obtain,

$$\pi_{\mathcal{R}}\boldsymbol{\varphi}_j = \sum_{\ell=1}^{6} (\boldsymbol{W}_{\mathcal{R}})_{(j\ell)} \left(\sum_{k=1}^{3n} \left[\boldsymbol{\chi}\left(\boldsymbol{r}_\ell\right)\right]_{(k)} \boldsymbol{\varphi}_k\right) = \sum_{\ell=1}^{6}\sum_{k=1}^{3n} (\boldsymbol{W}_{\mathcal{R}})_{(j\ell)} \left[\boldsymbol{\chi}\left(\boldsymbol{r}_\ell\right)\right]_{(k)} \boldsymbol{\varphi}_k \tag{57}$$

Let, $\boldsymbol{N}_{\mathcal{R}} \in \Re^{3n\times 6}$ be the matrix whose $(k,\ell)$th entry is $\left[\boldsymbol{\chi}\left(\boldsymbol{r}_\ell\right)\right]_{(k)}$, then

$$\pi_{\mathcal{R}}\boldsymbol{\varphi}_j = \sum_{k=1}^{3n} \left(\boldsymbol{N}_{\mathcal{R}}\boldsymbol{W}_{\mathcal{R}}^\top\right)_{(kj)} \boldsymbol{\varphi}_k \tag{58}$$

Finally, by comparing (58)with (54), we conclude that,

$$\boldsymbol{P}_{\mathcal{R}} = \boldsymbol{N}_{\mathcal{R}}\boldsymbol{W}_{\mathcal{R}}^\top \tag{59}$$

A similar derivation shows that the projection matrix $\boldsymbol{P}_{\mathcal{C}}$ can be expressed as:

$$\boldsymbol{P}_{\mathcal{C}} = \boldsymbol{N}_{\mathcal{C}}\boldsymbol{W}_{\mathcal{C}}^\top \tag{60}$$

where, $\boldsymbol{N}_{\mathcal{C}}, \boldsymbol{W}_{\mathcal{C}} \in \Re^{3n\times 6}$ with $(\boldsymbol{N}_{\mathcal{C}})_{(k\ell)} = \left[\boldsymbol{\chi}\left(\boldsymbol{c}_\ell\right)\right]_{(k)}$ and the $j$th row of $\boldsymbol{W}_{\mathcal{C}}$ given by

$$\left[\left\langle \boldsymbol{\epsilon}(\boldsymbol{\varphi}_j)\right\rangle_{(11)}, \left\langle \boldsymbol{\epsilon}(\boldsymbol{\varphi}_j)\right\rangle_{(22)}, \left\langle \boldsymbol{\epsilon}(\boldsymbol{\varphi}_j)\right\rangle_{(33)}, \left\langle \boldsymbol{\epsilon}(\boldsymbol{\varphi}_j)\right\rangle_{(12)}, \left\langle \boldsymbol{\epsilon}(\boldsymbol{\varphi}_j)\right\rangle_{(23)}, \left\langle \boldsymbol{\epsilon}(\boldsymbol{\varphi}_j)\right\rangle_{(31)}\right] \tag{61}$$

Note that the matrices $\boldsymbol{P}_{\mathcal{R}}$, $\boldsymbol{P}_{\mathcal{C}}$ are projections onto the range of $\boldsymbol{N}_{\mathcal{R}}$, $\boldsymbol{N}_{\mathcal{C}}$, respectively. We can also verify that $\boldsymbol{P}_{\mathcal{R}}\boldsymbol{P}_{\mathcal{C}} = \boldsymbol{P}_{\mathcal{R}}\boldsymbol{P}_{\mathcal{C}} = \boldsymbol{0}$. Using (27) and (53) we can find matrix representation of the projection $\pi_{\mathcal{P}}$ as

$$\boldsymbol{P}_{\mathcal{P}} = \boldsymbol{P}_{\mathcal{R}} + \boldsymbol{P}_{\mathcal{C}} \tag{62}$$

Finally, let us note that the null space of $\boldsymbol{P}_{\mathcal{P}}$ corresponds to the nodal representations of elements of $\mathcal{H}$.

### *4.2. Expressions for the stiffness matrix*

Using the projection matrices $\boldsymbol{P}_{\mathcal{R}}, \boldsymbol{P}_{\mathcal{C}}$ defined in previous section, we next obtain explicit expressions for the stiffness matrix, $\boldsymbol{K}_h^E$, associated with (41). We have

$$\left(\boldsymbol{K}_h^E\right)_{(jk)} = a_h^E(\boldsymbol{\varphi}_j, \boldsymbol{\varphi}_k) = a^E(\pi_{\mathcal{C}}\boldsymbol{\varphi}_j, \pi_{\mathcal{C}}\boldsymbol{\varphi}_k) + s^E(\boldsymbol{\varphi}_j - \pi_{\mathcal{P}}\boldsymbol{\varphi}_j, \boldsymbol{\varphi}_k - \pi_{\mathcal{P}}\boldsymbol{\varphi}_k) \tag{63}$$

To simplify the first term, let us first define the $6 \times 6$ matrix $\boldsymbol{D}$ whose entries are the normalized strain energies associated with uniform deformations, i.e.,

$$(\boldsymbol{D})_{(\ell m)} = \frac{1}{|E|} a^E(\boldsymbol{c}_\ell, \boldsymbol{c}_m), \qquad \ell, m = 1, \ldots, 6 \tag{64}$$

Using (22), we find that

$$\boldsymbol{D} = \begin{bmatrix} \boldsymbol{C}_{(1111)} & \boldsymbol{C}_{(1122)} & \boldsymbol{C}_{(1133)} & 2\boldsymbol{C}_{(1112)} & 2\boldsymbol{C}_{(1123)} & 2\boldsymbol{C}_{(1131)} \\ & \boldsymbol{C}_{(2222)} & \boldsymbol{C}_{(2233)} & 2\boldsymbol{C}_{(2212)} & 2\boldsymbol{C}_{(2223)} & 2\boldsymbol{C}_{(2231)} \\ & & \boldsymbol{C}_{(3333)} & 2\boldsymbol{C}_{(3312)} & 2\boldsymbol{C}_{(3323)} & 2\boldsymbol{C}_{(3331)} \\ & & & 4\boldsymbol{C}_{(1212)} & 4\boldsymbol{C}_{(1223)} & 4\boldsymbol{C}_{(1231)} \\ & \text{symm.} & & & 4\boldsymbol{C}_{(2323)} & 4\boldsymbol{C}_{(2331)} \\ & & & & & 4\boldsymbol{C}_{(3131)} \end{bmatrix} \tag{65}$$

which shows that $\boldsymbol{D}$ is only a function of the elasticity tensor $\boldsymbol{C}$ and does not depend on the geometry of the element $E$. For an isotropic material with Young's modulus $E_{\mathsf{Y}}$ and Poisson's ratio $\nu$, matrix $\boldsymbol{D}$ is given by:

$$\boldsymbol{D} = \frac{E_{\mathsf{Y}}}{(1+\nu)(1-2\nu)} \begin{bmatrix} 1-\nu & \nu & \nu & 0 & 0 & 0 \\ \nu & 1-\nu & \nu & 0 & 0 & 0 \\ \nu & \nu & 1-\nu & 0 & 0 & 0 \\ 0 & 0 & 0 & 2(1-2\nu) & 0 & 0 \\ 0 & 0 & 0 & 0 & 2(1-2\nu) & 0 \\ 0 & 0 & 0 & 0 & 0 & 2(1-2\nu) \end{bmatrix} \tag{66}$$

Using the expansion for $\pi_{\mathcal{C}} \boldsymbol{\varphi}_j$ similar to (55), we obtain the expression for the first term of the stiffness matrix as

$$\begin{aligned} a^E(\pi_{\mathcal{C}}\boldsymbol{\varphi}_j, \pi_{\mathcal{C}}\boldsymbol{\varphi}_k) &= a^E\Big(\sum_{\ell=1}^{6} (\boldsymbol{W}_{\mathcal{C}})_{(j\ell)}\, \boldsymbol{c}_\ell, \sum_{m=1}^{6} (\boldsymbol{W}_{\mathcal{C}})_{(km)}\, \boldsymbol{c}_m\Big) \\ &= \sum_{\ell=1}^{6}\sum_{\ell=1}^{6} (\boldsymbol{W}_{\mathcal{C}})_{(j\ell)} \left[a^E(\boldsymbol{c}_\ell, \boldsymbol{c}_m)\right] (\boldsymbol{W}_{\mathcal{C}})_{(km)} \\ &= |E| \left(\boldsymbol{W}_{\mathcal{C}} \boldsymbol{D} \boldsymbol{W}_{\mathcal{C}}^\top\right)_{(jk)} \end{aligned} \tag{67}$$

As for the second term, we first define the $3n \times 3n$ matrix $\boldsymbol{S}^E$ whose $(j,k)$th entry is $s^E(\boldsymbol{\varphi}_j, \boldsymbol{\varphi}_k)$. For example, for the choice of $s^E$ given by (47), we have $\boldsymbol{S}^E = \alpha^E \boldsymbol{I}_{3n}$. Noting that

$$\boldsymbol{\varphi}_j - \pi_{\mathcal{P}}\boldsymbol{\varphi}_j = \sum_{k=1}^{3n} (\boldsymbol{I} - \boldsymbol{P}_{\mathcal{P}})_{(kj)}\, \boldsymbol{\varphi}_k \tag{68}$$

we can write

$$s^E(\boldsymbol{\varphi}_j - \pi_{\mathcal{P}}\boldsymbol{\varphi}_j, \boldsymbol{\varphi}_k - \pi_{\mathcal{P}}\boldsymbol{\varphi}_k) = \left[(\boldsymbol{I} - \boldsymbol{P}_{\mathcal{P}})^\top \boldsymbol{S}^E (\boldsymbol{I} - \boldsymbol{P}_{\mathcal{P}})\right]_{(jk)} \tag{69}$$

and so the stiffness matrix is given by

$$\boldsymbol{K}_h^E = |E|\, \boldsymbol{W}_{\mathcal{C}} \boldsymbol{D} \boldsymbol{W}_{\mathcal{C}}^\top + (\boldsymbol{I} - \boldsymbol{P}_{\mathcal{P}})^\top \boldsymbol{S}^E (\boldsymbol{I} - \boldsymbol{P}_{\mathcal{P}}) \tag{70}$$

As seen from the above expression, the task of computing the stiffness matrix reduces to computing four matrices $\boldsymbol{N}_{\mathcal{R}}$, $\boldsymbol{N}_{\mathcal{C}}$, $\boldsymbol{W}_{\mathcal{R}}$, and $\boldsymbol{W}_{\mathcal{C}}$. We next further break down these calculation and comment on the computational effort needed for the method.

### 4.3. Calculation of matrices $\boldsymbol{N}_{\mathcal{R}}$, $\boldsymbol{N}_{\mathcal{C}}$, $\boldsymbol{W}_{\mathcal{R}}$, and $\boldsymbol{W}_{\mathcal{C}}$

Lets first concentrate on the matrices $\boldsymbol{N}_{\mathcal{R}}$ and $\boldsymbol{N}_{\mathcal{C}}$ which are essentially the nodal representations of the bases of $\mathcal{R}$ and $\mathcal{C}$, respectively. Referring to (21), we can see that, for $i = 1, \ldots, n$, the block of $3i-2$ to $3i$ rows of $\boldsymbol{N}_{\mathcal{R}}$, which are associated with the $i$th vertex of the element, is given by:

$$\begin{bmatrix} 1 & 0 & 0 & (\boldsymbol{x}_i - \overline{\boldsymbol{x}})_{(2)} & 0 & -(\boldsymbol{x}_i - \overline{\boldsymbol{x}})_{(3)} \\ 0 & 1 & 0 & -(\boldsymbol{x}_i - \overline{\boldsymbol{x}})_{(1)} & (\boldsymbol{x}_i - \overline{\boldsymbol{x}})_{(3)} & 0 \\ 0 & 0 & 1 & 0 & -(\boldsymbol{x}_i - \overline{\boldsymbol{x}})_{(2)} & (\boldsymbol{x}_i - \overline{\boldsymbol{x}})_{(1)} \end{bmatrix} \tag{71}$$

Similarly, the block $3i-2$ to $3i$ rows of $\boldsymbol{N}_{\mathcal{C}}$ can be expressed as:

$$\begin{bmatrix} (\boldsymbol{x}_i - \overline{\boldsymbol{x}})_{(1)} & 0 & 0 & (\boldsymbol{x}_i - \overline{\boldsymbol{x}})_{(2)} & 0 & (\boldsymbol{x}_i - \overline{\boldsymbol{x}})_{(3)} \\ 0 & (\boldsymbol{x}_i - \overline{\boldsymbol{x}})_{(2)} & 0 & (\boldsymbol{x}_i - \overline{\boldsymbol{x}})_{(1)} & (\boldsymbol{x}_i - \overline{\boldsymbol{x}})_{(3)} & 0 \\ 0 & 0 & (\boldsymbol{x}_i - \overline{\boldsymbol{x}})_{(3)} & 0 & (\boldsymbol{x}_i - \overline{\boldsymbol{x}})_{(2)} & (\boldsymbol{x}_i - \overline{\boldsymbol{x}})_{(1)} \end{bmatrix} \tag{72}$$

To facilitate the description of the block of $\boldsymbol{W}_{\mathcal{R}}$ and $\boldsymbol{W}_{\mathcal{C}}$ associated with the $i$th vertex, we first define the vector $\boldsymbol{q}_i$ as:

$$\boldsymbol{q}_i = \frac{1}{2|E|} \int_{\partial E} \varphi_i \boldsymbol{n} \mathrm{d}\boldsymbol{s} \tag{73}$$

Recall that $\varphi_i$ is the scalar barycentric coordinate based on which the element basis functions are defined (cf., (9)). Since $\varphi_i$ vanishes on $\mathcal{F}_i^c$, we have

$$\boldsymbol{q}_i = \frac{1}{2|E|} \sum_{F \in \mathcal{F}_i} \left( \int_F \varphi_i \mathrm{d}\boldsymbol{s} \right) \boldsymbol{n}_{F,E} \tag{74}$$

where $\boldsymbol{n}_{F,E}$ is the normal to the face $F$ pointing outwards with respect to element $E$. Using (9), the definition of $\boldsymbol{W}_{\mathcal{R}}$ in (56), and the identity (31), we find the block of $3i-2$ to $3i$ rows of $\boldsymbol{W}_{\mathcal{R}}$ to be

$$\begin{bmatrix} 1/n & 0 & 0 & (\boldsymbol{q}_i)_{(2)} & 0 & -(\boldsymbol{q}_i)_{(3)} \\ 0 & 1/n & 0 & -(\boldsymbol{q}_i)_{(1)} & (\boldsymbol{q}_i)_{(3)} & 0 \\ 0 & 0 & 1/n & 0 & -(\boldsymbol{q}_i)_{(2)} & (\boldsymbol{q}_i)_{(1)} \end{bmatrix} \tag{75}$$

A similar approach can be used to identify the block of $3i-2$ to $3i$ rows of $\boldsymbol{W}_{\mathcal{C}}$ as

$$\begin{bmatrix} 2(\boldsymbol{q}_i)_{(1)} & 0 & 0 & (\boldsymbol{q}_i)_{(2)} & 0 & (\boldsymbol{q}_i)_{(3)} \\ 0 & 2(\boldsymbol{q}_i)_{(2)} & 0 & (\boldsymbol{q}_i)_{(1)} & (\boldsymbol{q}_i)_{(3)} & 0 \\ 0 & 0 & 2(\boldsymbol{q}_i)_{(3)} & 0 & (\boldsymbol{q}_i)_{(2)} & (\boldsymbol{q}_i)_{(1)} \end{bmatrix} \tag{76}$$

The calculation of the matrices $\boldsymbol{W}_{\mathcal{R}}$ and $\boldsymbol{W}_{\mathcal{C}}$ thus boils down to computation of the surface integrals $\int_F \varphi_i \mathrm{d}\boldsymbol{s}$ as well the unit normal vectors for each face in the mesh. Two additional data structures for $\mathcal{T}_h$, beyond the usual vertex list and vertex-element connectivity matrix, are needed to facilitate these calculations. The first contains the list of faces in the mesh along with vertices incident on each face. The orientation of each face $F$, and subsequently its unit normal vector $\boldsymbol{n}_F$, is chosen and fixed once and for all. The second data structure contains the element-face connectivity information as well as the orientation of the outer normal $\boldsymbol{n}_{F,E}$ associated with element $E$. Note that either $\boldsymbol{n}_{F,E} = \boldsymbol{n}_F$ or $\boldsymbol{n}_{F,E} = -\boldsymbol{n}_F$, depending on the initial choice of $\boldsymbol{n}_F$. We should point out that the construction of these data structures directly from the standard vertex-element connectivity and vertex list, is possible since the faces of an element can be identified by inspecting which groups of vertices that lie on the same plane. If the mesh consists of convex polyhedra, the process can be greatly simplified. For example, each internal faces of such a mesh is shared exactly between two elements. This topological insight can be used to identify the internal faces through an inspection of the vertex-element connectivity matrix.

Regarding the surface integrals $\int_F \varphi_i \mathrm{d}\boldsymbol{s}$, if boundary barycentric coordinates of [3] are used in the construction of $\mathcal{W}(E)$, the value of this surface integral can be computed exactly from the geometry of $F$ (see the appendix). In general, we can use the nodal quadrature scheme presented in Section 3.2 to compute a first-order approximation. Note that the nodal quadrature does not distinguish between the choice of basis functions and thus the value of the surface is purely a function of the geometry of $F$[7]. Such an approximation in fact amounts to using the following projection maps

$$\pi_{\mathcal{R}}\boldsymbol{v} = \overline{\boldsymbol{v}} + \left[\frac{1}{2\,|E|}\sum_{F\subseteq\partial E} \fint_F (\boldsymbol{v}\otimes\boldsymbol{n} - \boldsymbol{n}\otimes\boldsymbol{v})\,\mathrm{d}\boldsymbol{s}\right](\boldsymbol{x}-\overline{\boldsymbol{x}}) \tag{77}$$

$$\pi_{\mathcal{C}}\boldsymbol{v} = \left[\frac{1}{2\,|E|}\sum_{F\subseteq\partial E} \fint_F (\boldsymbol{v}\otimes\boldsymbol{n} + \boldsymbol{n}\otimes\boldsymbol{v})\,\mathrm{d}\boldsymbol{s}\right](\boldsymbol{x}-\overline{\boldsymbol{x}}) \tag{78}$$

in (41). Observe that it is necessary for the quadrature rule to exactly integrate linear polynomials in order for these maps to respect (23)-(26). We can show that the consistency error in the energy bilinear form (i.e., the first term in (15)) will remain $\mathcal{O}(h)$. While the orthogonality conditions (36)-(38), and subsequently (43), are satisfied only asymptotically, we will show in Section 5.1 that the global patch are nevertheless will be passed exactly.

### *4.4. Structure the stiffness matrix*

It is insightful to examine the effects of strain energy approximation in VEM from an algebraic point of view by considering the structure of the resulting stiffness matrix. Such a perspective is fundamental to the development of corresponding MFD formulations that do not directly utilize the existence of underlying basis functions (see, for example, [18, 7]).

To this effect, let us again consider an element $E$ and denote its *exact* stiffness matrix by $\boldsymbol{K}^E$, i.e., $\boldsymbol{K}^E_{(jk)} = a^E(\boldsymbol{\varphi}_j, \boldsymbol{\varphi}_k)$. We assume that the *diagonal* stability term (47) is used in the definition of the discrete bilinear form and the VEM stiffness $\boldsymbol{K}^E_h$. Complementing the bases for space of rigid body motions and constant strains, we will next select a basis $\boldsymbol{h}_1, \ldots, \boldsymbol{h}_{3n-12}$ for the space of higher-order modes[8] $\mathcal{H} = \mathcal{H}(E)$. This basis will be chosen such that the corresponding matrix $\boldsymbol{N}_{\mathcal{H}} \in \Re^{3n\times(3n-12)}$, i.e., the matrix whose $\ell$th column is $\boldsymbol{\chi}(\boldsymbol{h}_\ell)$, satisfies the following properties:

$$\boldsymbol{N}_{\mathcal{H}}^\top \boldsymbol{N}_{\mathcal{H}} = \boldsymbol{I}_{3n-12}, \qquad \boldsymbol{N}_{\mathcal{H}}^\top \boldsymbol{K}^E \boldsymbol{N}_{\mathcal{H}} = \mathrm{diag}(\lambda_1^E, \ldots, \lambda_{3n-12}^E) \tag{79}$$

[7] In fact, $\fint_F \varphi_i \mathrm{d}\boldsymbol{s}$ is the simply the area of the quadrilateral formed by the associated vertex $\boldsymbol{x}_i$, the midpoint of edges incident on $\boldsymbol{x}_i$ and the centroid of $F$.

[8] Here, we assume that the exact projections maps are available andused to define $\mathcal{H}$.

Observe that $\lambda_\ell^E$ is the strain energy associated with $\boldsymbol{h}_\ell$, i.e., $a^E(\boldsymbol{h}_\ell, \boldsymbol{h}_\ell) = \lambda_\ell^E$. To see how such a basis can be selected, we start with matrix $\tilde{\boldsymbol{N}}_{\mathcal{H}} \in \Re^{3n\times(3n-12)}$ whose columns form an orthonormal basis for the null space of $\boldsymbol{P}_{\mathcal{P}}$. We choose $\boldsymbol{U}$ to be the orthogonal matrix whose columns are eigenvectors of $\tilde{\boldsymbol{N}}_{\mathcal{H}}^\top \boldsymbol{K}^E \tilde{\boldsymbol{N}}_{\mathcal{H}}$ and set $\boldsymbol{N}_{\mathcal{H}} \doteq \tilde{\boldsymbol{N}}_{\mathcal{H}} \boldsymbol{U}$. Since $\boldsymbol{P}_{\mathcal{P}} \boldsymbol{N}_{\mathcal{H}} = \boldsymbol{P}_{\mathcal{P}} \tilde{\boldsymbol{N}}_{\mathcal{H}} \boldsymbol{U} = \boldsymbol{0}$, the columns of $\boldsymbol{N}_{\mathcal{H}}$ indeed corresponds to deformation states in $\mathcal{H}$. One can verify that conditions in (79) hold in this case.

Next, we define a $3n \times 3n$ matrix given by

$$\boldsymbol{N} = [\boldsymbol{N}_{\mathcal{R}}, \boldsymbol{N}_{\mathcal{C}}, \boldsymbol{N}_{\mathcal{H}}] \doteq [\boldsymbol{\chi}(\boldsymbol{r}_1), \ldots, \boldsymbol{\chi}(\boldsymbol{r}_6), \boldsymbol{\chi}(\boldsymbol{c}_1), \ldots, \boldsymbol{\chi}(\boldsymbol{r}_6), \boldsymbol{\chi}(\boldsymbol{h}_1), \ldots, \boldsymbol{\chi}(\boldsymbol{h}_{3n-12})] \tag{80}$$

which is invertible since the bases for rigid body motion, uniform and higher-order deformations are linearly independent. Owing to the energy-orthogonality of $\mathcal{C}$ and $\mathcal{H}$, the exact stiffness matrix form has a block diagonal structure with respect to basis defined by the columns of $\boldsymbol{N}$. More specifically, we have

$$\begin{aligned} \boldsymbol{N}^\top \boldsymbol{K}^E \boldsymbol{N} &= \begin{bmatrix} \boldsymbol{N}_{\mathcal{R}}^\top \boldsymbol{K}^E \boldsymbol{N}_{\mathcal{R}} & \boldsymbol{N}_{\mathcal{R}}^\top \boldsymbol{K}^E \boldsymbol{N}_{\mathcal{C}} & \boldsymbol{N}_{\mathcal{R}}^\top \boldsymbol{K}^E \boldsymbol{N}_{\mathcal{H}} \\ \boldsymbol{N}_{\mathcal{C}}^\top \boldsymbol{K}^E \boldsymbol{N}_{\mathcal{R}} & \boldsymbol{N}_{\mathcal{C}}^\top \boldsymbol{K}^E \boldsymbol{N}_{\mathcal{C}} & \boldsymbol{N}_{\mathcal{C}}^\top \boldsymbol{K}^E \boldsymbol{N}_{\mathcal{H}} \\ \boldsymbol{N}_{\mathcal{H}}^\top \boldsymbol{K}^E \boldsymbol{N}_{\mathcal{R}} & \boldsymbol{N}_{\mathcal{H}}^\top \boldsymbol{K}^E \boldsymbol{N}_{\mathcal{C}} & \boldsymbol{N}_{\mathcal{H}}^\top \boldsymbol{K}^E \boldsymbol{N}_{\mathcal{H}} \end{bmatrix} \\ &= \begin{bmatrix} \boldsymbol{0} & \boldsymbol{0} & \boldsymbol{0} \\ \boldsymbol{0} & |E|\,\boldsymbol{D} & \boldsymbol{0} \\ \boldsymbol{0} & \boldsymbol{0} & \mathrm{diag}(\lambda_1^E, \ldots, \lambda_{3n-12}^E) \end{bmatrix} \end{aligned} \tag{81}$$

Similarly, one obtains the following transformation of the VEM stiffness matrix under the same change of basis:

$$\boldsymbol{N}^\top \boldsymbol{K}_h^E \boldsymbol{N} = \begin{bmatrix} \boldsymbol{0} & \boldsymbol{0} & \boldsymbol{0} \\ \boldsymbol{0} & |E|\,\boldsymbol{D} & \boldsymbol{0} \\ \boldsymbol{0} & \boldsymbol{0} & \alpha^E \boldsymbol{I}_{3n-12} \end{bmatrix} \tag{82}$$

Compared to (81), we can see that the only difference lies in the strain energy assigned to the higher-order modes $\boldsymbol{h}_1, \ldots, \boldsymbol{h}_{3n-12}$. With the choice of stabilizing bilinear form (47), the energy of these modes are assumed to identically equal to $\alpha^E$.

This observation leads us to the question of selection of an appropriate value for the parameter $\alpha^E$. While there are other possible approaches to addressing this question, we proceed as follows. Because $s^E$ can be viewed an approximation to the exact strain energy $a^E$, we can compare the energy it assigns to uniform deformation with their exact energy. Observing that

$$s^E(\boldsymbol{c}_\ell, \boldsymbol{c}_m) = \left(\boldsymbol{N}_{\mathcal{C}}^\top \boldsymbol{S}^E \boldsymbol{N}_{\mathcal{C}}\right)_{(\ell m)} = \alpha^E \left(\boldsymbol{N}_{\mathcal{C}}^\top \boldsymbol{N}_{\mathcal{C}}\right)_{(\ell m)} \tag{83}$$

and $a^E(\boldsymbol{c}_\ell, \boldsymbol{c}_m) = |E|\,\boldsymbol{D}_{(\ell m)}$, the scaling coefficient can be chosen such that $\alpha^E \boldsymbol{N}_{\mathcal{C}}^\top \boldsymbol{N}_{\mathcal{C}}$ and $|E|\,\boldsymbol{D}$ are comparable. Equating the trace of these matrix suggests the following relation for the scaling parameter

$$\alpha^E = \gamma \alpha_\star^E \tag{84}$$

where

$$\alpha_\star^E \doteq \frac{|E|\,\mathrm{trace}\,(\boldsymbol{D})}{\mathrm{trace}\left(\boldsymbol{N}_\mathcal{C}^\top \boldsymbol{N}_\mathcal{C}\right)} \tag{85}$$

The coefficient $\gamma$, taken to be independent of $E$, is introduced here in order to facilitate the study of the influence of $\alpha^E$ on the accuracy of the numerical solutions. As expected, this expression depends on the material properties, through the trace of $\boldsymbol{D}$, and the geometry of $E$ through the appearance of $\boldsymbol{N}_\mathcal{C}$.

The above argument indicates that a reasonable value for $\gamma$ should be close to one. In order to validate this assertion, we have computed the strain energy associated with basis for $\mathcal{H}$, i.e., $\lambda_1^E, \ldots, \lambda_{3n-12}^E$, of a few representative polyhedra shown in Fig. 5. The exact stiffness matrix here is computed, by means of very high order quadrature scheme, for an element space $\mathcal{W}(E)$ defined using maximum entropy basis functions[9]. The spectrum is normalized by $\alpha_\star^E$ in order facilitate the comparison with coefficient $\gamma$. We can see that the spectrum (c.f. Fig. 6) can be relatively broad especially for the irregularly-shaped distorted polyhedra with small edges. However, the average value of the normalized energies is $\mathcal{O}(1)$ for the element geometries considered here. The histogram in Fig. 7 shows the distribution of this average value for the elements in a centroidal Voronoi mesh with 100 polyhedral elements. We can see a significant clustering around 1 confirming that (84) with $\gamma = 1$ is a reasonable choice for the strain energy of higher order modes.

## 5. Numerical studies

In this section, we evaluate the performance of polyhedral discretizations using VEM through several numerical studies. The accuracy and convergence of the numerical solutions are assessed using two measures of error in computed displacement and stress fields. The relative error in displacements is defined using the volumetric nodal quadrature rule of Section 3 as:

$$\mathsf{e}_{\boldsymbol{u}} \doteq \left[\frac{\sum_{E\in\mathcal{T}_h} \oint_E \left|\boldsymbol{u} - \tilde{\boldsymbol{u}}_h\right|^2 \mathrm{d}\boldsymbol{x}}{\sum_{E\in\mathcal{T}_h} \oint_E \left|\boldsymbol{u}\right|^2 \mathrm{d}\boldsymbol{x}}\right]^{1/2} \tag{86}$$

Recall that $\tilde{\boldsymbol{u}}_h$ is the VEM solution and $\boldsymbol{u}$ is the exact solution, assumed to be sufficiently smooth for the quadrature to make sense. Note that $\mathsf{e}_{\boldsymbol{u}}$ serves as an approximation to the $L^2$-norm error $\|\boldsymbol{u} - \tilde{\boldsymbol{u}}_h\|_{0,\Omega}$, which cannot be computed without access to the basis functions.

Similarly, we define a discrete error measure for the stress field since the *raw field* $\boldsymbol{\sigma}(\tilde{\boldsymbol{u}}_h)$ are not readily available. Motivating by the fact that $\boldsymbol{\epsilon}(\pi_{\mathcal{C}(E)}\boldsymbol{v})$ is the volume average of $\boldsymbol{\epsilon}(\boldsymbol{v})$ over element $E$ (cf. (33)) and therefore its best constant approximation, we define a element-wise constant stress field on $\mathcal{T}_h$, denoted by $\boldsymbol{\sigma}_h(\boldsymbol{v})$, such that[10]

$$\boldsymbol{\sigma}_h(\boldsymbol{v})\Big|_E = \frac{1}{|E|} \int_E \boldsymbol{\sigma}(\boldsymbol{v}) \mathrm{d}\boldsymbol{x} = \boldsymbol{\sigma}(\pi_{\mathcal{C}} \boldsymbol{v}) \tag{87}$$

[9] In the case of cube, the underlying space $\mathcal{W}$ coincides with the usual finite element space of trilinear functions.

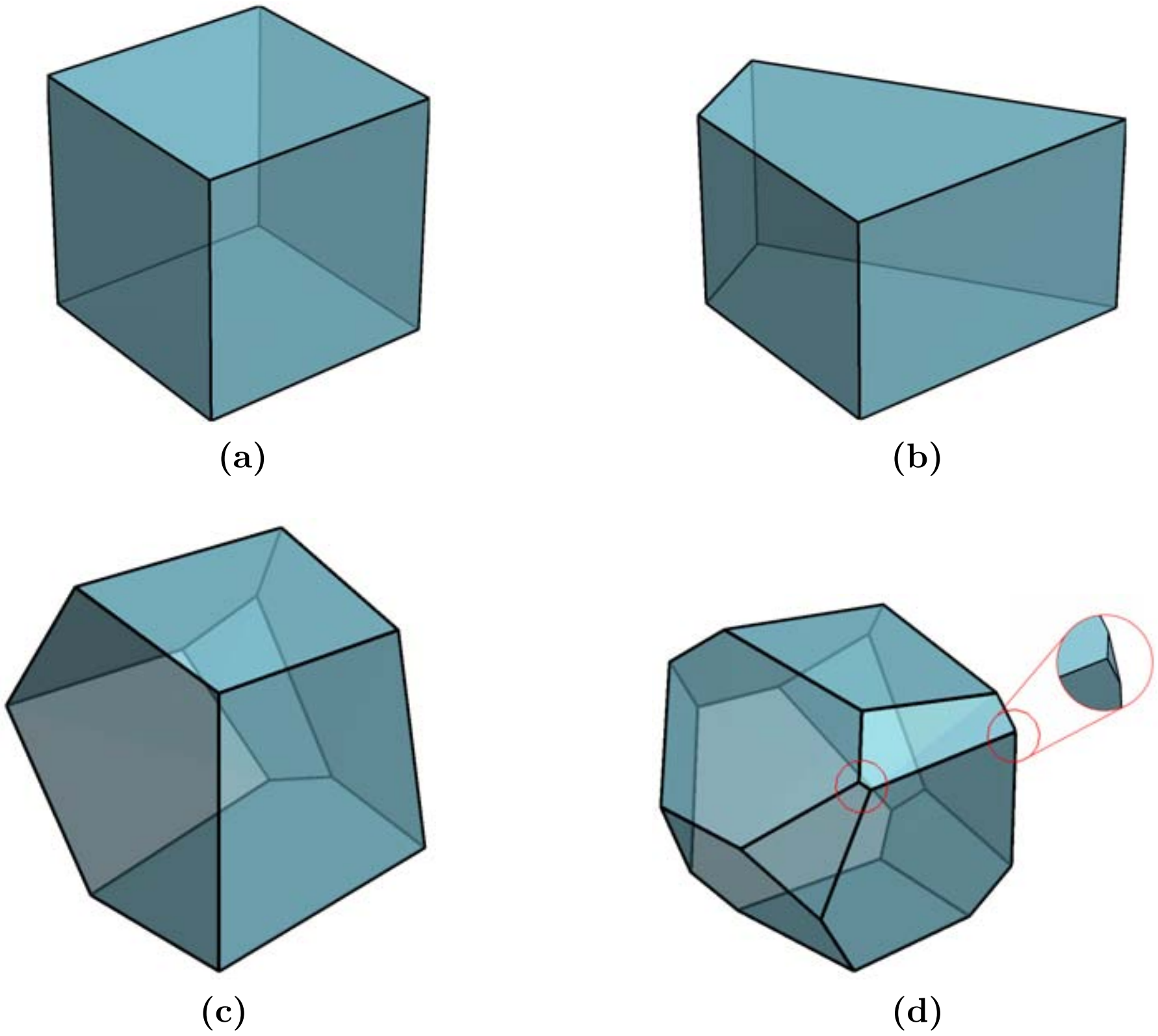

**Figure 5:** Representative polyhedra for the study of energy of higher-order modes $\lambda_\ell^E$ (a) Uniform hexahedron. (b) Distorted hexahedron. (c) Polyhedron with 12 vertices (d) Polyhedron with 24 vertices. Red circles highlight small edges.

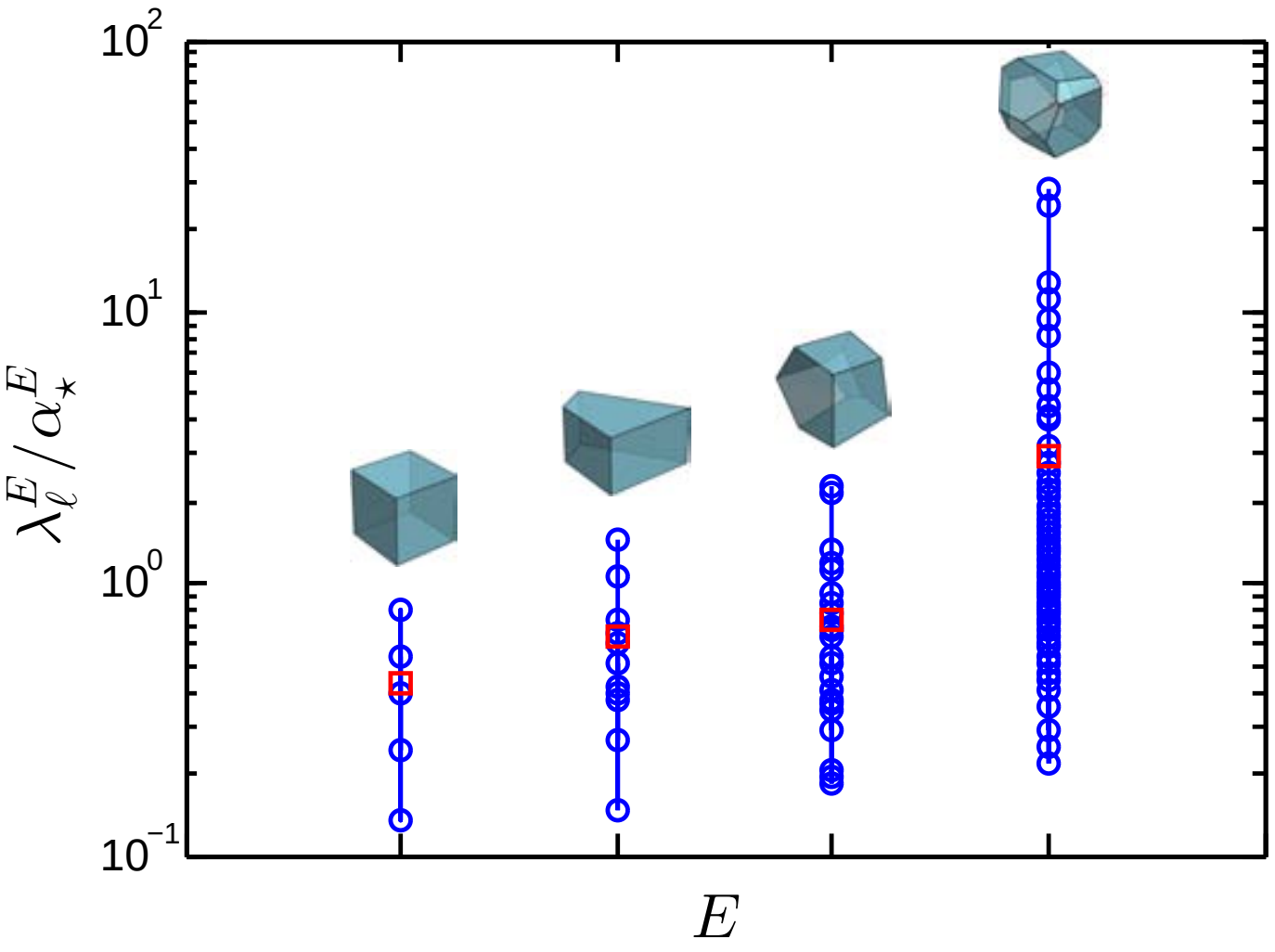

**Figure 6:** Distribution of normalized strain energies associated with higher-order modes, $\lambda_\ell^E/\alpha_\star^E$, for the polyhedra shown in Fig. 5. The red square box within each spectrum indicates the average value.

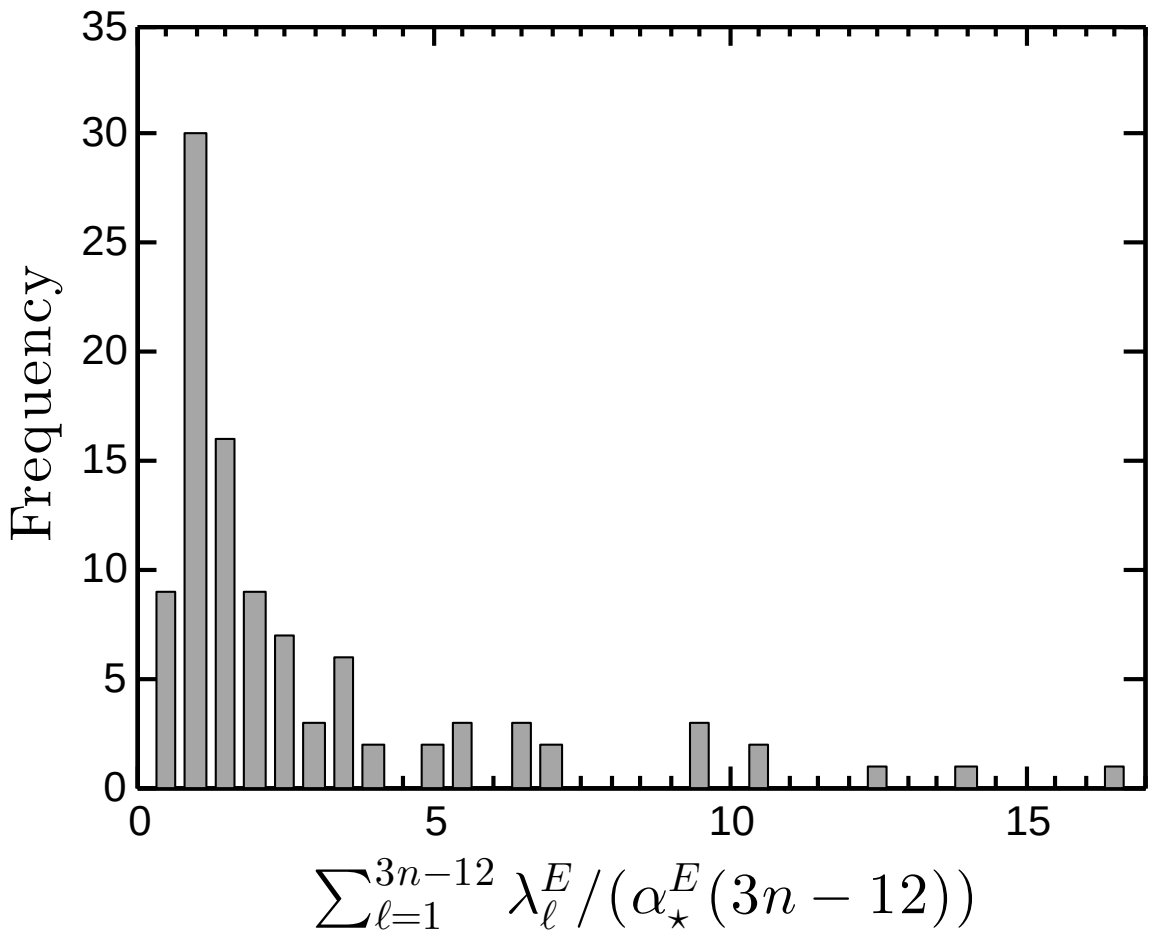

**Figure 7:** Frequency distribution of the average normalized strain energy $\sum_{\ell=1}^{3n-12} \lambda_\ell^E / (\alpha_\star^E (3n-12))$ for a Centroidal Voronoi mesh consisting of 100 polyhedral elements.

This will serve as a surrogate to $\boldsymbol{\sigma}(\boldsymbol{v})$ and, accordingly, the following measure of error is considered

$$\mathsf{e}_{\boldsymbol{\sigma}} \doteq \frac{\|\boldsymbol{\sigma}(\boldsymbol{u}) - \boldsymbol{\sigma}_h(\tilde{\boldsymbol{u}}_h)\|_{0,\Omega}}{\|\boldsymbol{\sigma}(\boldsymbol{u})\|_{0,\Omega}} \tag{88}$$

We evaluate these integrals numerically with a high-order quadrature obtained as follows: each polyhedral element is divided into pyramids and a forth-order Gauss rule consisting of 64 integration points is used over each pyramid (cf. [34]).

We use the open source MATLAB toolbox, Multi-Parametric Toolbox (MPT) [26], for generating the polyhedral meshes. Two types of polyhedral meshes based on Voronoi tessellations are used in our study: random Voronoi meshes (abbreviated by RND), which are formed by a random set of generating seeds; and more uniform centroidal Voronoi tessellations (abbreviated by CVT) that are obtained from a set of seeds that coincide with centroids of the resulting Voronoi cells. The CVT meshes are generated using Lloyd's algorithm following the approach outlined in [38]. Both types of meshes consist only of convex polyhedra. Finally, in all the numerical presented here, the local discrete bilinear forms are based on choice of $s^E$ in (47) with $\alpha^E$ given by (84).

### *5.1. Displacement patch test*

We start with the displacement patch test on the unit cube $\Omega = (0,1)^3$ discretized using CVT and RND meshes with different number of polyhedrons. We have also considered

[10] We can express this average element stress as $\boldsymbol{\sigma}_h(\boldsymbol{v})|_E = \sum_{\ell=1}^{6} [\boldsymbol{W}_{\mathcal{C}}^\top \boldsymbol{\chi}(\boldsymbol{v})]_\ell \boldsymbol{\sigma}(\boldsymbol{c}_\ell)$.

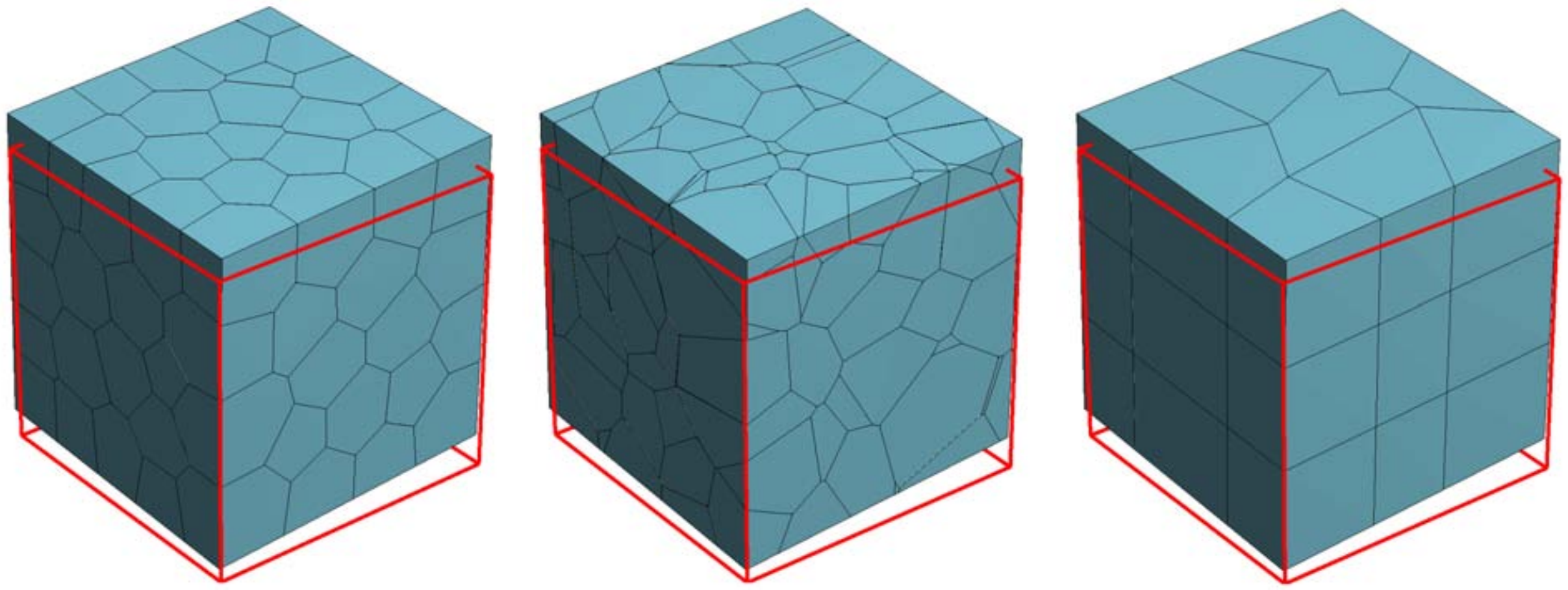

**(a)** **(b)** **(c)**

**Figure 8:** Illustration of the displacement patch test with exact solution $\boldsymbol{p} = [2x_1 + x_2 + 3x_3 + 1, 3x_1 + 4x_2 + 2x_3 + 2, 4x_1 + 3x_2 + x_3 + 3]^\top/100$ (a) CVT mesh, (b) RND mesh, and (c) polyhedral mesh containing non-convex elements. The red box represents $\Omega$ and the mesh illustrates the deformed configuration.

meshes containing non-convex elements. The exact solution for the patch test is an arbitrary linear displacement field $\boldsymbol{u} = \boldsymbol{p} \in \mathcal{P}(\Omega)$. We consider the case where $\Gamma_{\boldsymbol{t}} = \emptyset$ and $\boldsymbol{g} = \boldsymbol{p}|_{\partial\Omega}$ is applied to the entire boundary (note that the body forces are absent, i.e., $\boldsymbol{b} = \boldsymbol{0}$). Figure 8 shows the deformed configuration of the three representative meshes tested. The relative displacement and stress errors close to machine precision levels, for a large range of $\gamma$ values, are observed for CVT meshes indicating that VEM passes displacement patch test. For random Voronoi meshes the errors are approximately one order of magnitude higher than CVT mesh errors but still close to machine precision levels.

These results are consistent with the theoretical discussion thus far when the projection maps are computed exactly (e.g., when the element space is obtained using boundary coordinates of [3]). However, we have also observed the satisfaction of the global patch test even when nodal quadrature is used in the definition of the projection maps (cf. (77) and (78)). In this case, the element level condition of exactness of strain energy, i.e., (43), may not be satisfied exactly. Nevertheless we can prove directly that the global patch test should be passed. First note that, for $\boldsymbol{p} \in \mathcal{P}(E)$ and $\boldsymbol{v} \in \mathcal{W}(E)$, and using projection maps obtained from nodal quadrature, we have

$$\begin{aligned} a_h^E(\boldsymbol{p}, \boldsymbol{v}) &= a^E(\boldsymbol{p}, \pi_{\mathcal{C}}\boldsymbol{v}) \\ &= \boldsymbol{\sigma}(\boldsymbol{p}) : \left[\int_E \boldsymbol{\epsilon}(\pi_{\mathcal{C}}\boldsymbol{v})\mathrm{d}\boldsymbol{x}\right] \\ &= \boldsymbol{\sigma}(\boldsymbol{p}) : \left[\sum_{F\subseteq\partial E} \frac{1}{2}\int_F (\boldsymbol{v}\otimes\boldsymbol{n} + \boldsymbol{n}\otimes\boldsymbol{v})\,\mathrm{d}\boldsymbol{s}\right] \\ &= \sum_{F\subseteq\partial E}\int_F \boldsymbol{v}\cdot\boldsymbol{\sigma}(\boldsymbol{p})\boldsymbol{n}\mathrm{d}\boldsymbol{s} \qquad \text{(using symmetry of } \boldsymbol{\sigma}(\boldsymbol{p})\text{)} \end{aligned} \tag{89}$$

Now, let us consider a general patch test with exact solution $\boldsymbol{u} = \boldsymbol{p} \in \mathcal{P}(\Omega)$. Corresponding tractions $\boldsymbol{t} = \boldsymbol{\sigma}(\boldsymbol{p})\boldsymbol{n}$ are imposed on $\Gamma_{\boldsymbol{t}}$, and $\boldsymbol{g} = \boldsymbol{p}|_{\Gamma_{\boldsymbol{u}}}$ is applied to the remainder of the boundary. For an arbitrary test function $\boldsymbol{v} \in \mathcal{V}_h^{\boldsymbol{0}}$, we see that

$$a_h(\boldsymbol{p}, \boldsymbol{v}) = \sum_{E \in \mathcal{T}_h} a_h^E(\boldsymbol{p}, \boldsymbol{v}) = \sum_{E \in \mathcal{T}_h} \sum_{F \subseteq \partial E} \fint_F \boldsymbol{v} \cdot \boldsymbol{\sigma}(\boldsymbol{p})\boldsymbol{n} \mathrm{d}\boldsymbol{s} = \sum_{F \subseteq \Gamma_{\boldsymbol{t}}} \fint_F \boldsymbol{v} \cdot \boldsymbol{t} \mathrm{d}\boldsymbol{s} = f_h(\boldsymbol{v}) \tag{90}$$

**Table 1:** Comparison between maximum entropy finite element and VEM for displacement patch test. The second and third columns show the relative displacement errors $\mathsf{e}_{\boldsymbol{u}}$.

| Number of elements | Maximum entropy | VEM |
|---|---|---|
| 50 | $8.14 \times 10^{-2}$ | $5.72 \times 10^{-15}$ |
| 100 | $9.24 \times 10^{-2}$ | $1.91 \times 10^{-14}$ |
| 200 | $6.63 \times 10^{-2}$ | $2.66 \times 10^{-14}$ |

In the second to last equality, we have used the fact that boundary integrals on the internal faces cancel out and $\boldsymbol{v} = \boldsymbol{0}$ on $\Gamma_{\boldsymbol{u}}$. Since $\boldsymbol{p} \in \mathcal{V}_h^{\boldsymbol{g}}$, this shows that $\tilde{\boldsymbol{u}}_h = \boldsymbol{p}$ is the unique solution to the discrete problem and the global patch test is passed.

For the sake of comparison, we repeat the patch test using a finite element method based on the conforming maximum entropy basis functions [36] and numerical integration of the stiffness matrix. The errors in the patch test performed for a sequence of meshes is indicative of the consistency error introduced by the quadrature in the elemental and global strain energies. Numerical integration in each element is carried out over by first partitioning into tetrahedra[11] and then standard quadrature rules to each tetrahedron are applied. On a CVT mesh of 50 elements, the observed relative errors in displacement are $1.84 \times 10^{-1}$, $8.14 \times 10^{-2}$ and $4.93 \times 10^{-2}$ when first, second and forth order quadrature rules (consisting of 1, 4 and 11 points) are used for each tetrahedron. Table 1 shows the errors for a sequence of CVT meshes with the second order rule indicating that the errors persist under mesh refinement. As discussed in detail in [37], the failure to satisfy patch test (at least asymptotically under mesh refinement) can place a limit on the accuracy that can be achieved in general by the method.

### *5.2. Shear-loaded beam*

Next, we study the performance of VEM for the cantilever beam loaded in shear. The domain $\Omega$ for this problem is $(-1, 1) \times (-1, 1) \times (0, L)$, which is occupied by an isotropic material with Young's modulus $E_{\mathsf{Y}}$ and Poisson's ratio $\nu$, and is subjected to constant tractions, given by $\boldsymbol{t} = [0, -F, 0]^\top$, on the face passing through the origin. The expressions

[11] Each polyhedral element is divided into tetrahedra using the element center, face centers and vertex locations.

for stresses are available in Barber [4] and repeated here for completeness:

$$\begin{aligned}
&\boldsymbol{\sigma}_{(11)} = \boldsymbol{\sigma}_{(22)} = \boldsymbol{\sigma}_{(12)} = 0, \quad \boldsymbol{\sigma}_{(33)} = \frac{3F}{4}\boldsymbol{x}_{(2)}\boldsymbol{x}_{(3)} \\
&\boldsymbol{\sigma}_{(31)} = \frac{3F\nu}{2\pi^2(1+\nu)} \sum_{n=1}^{\infty} \frac{(-1)^n}{n^2\cosh(n\pi)} \sin\left(n\pi \boldsymbol{x}_{(1)}\right) \sinh\left(n\pi \boldsymbol{x}_{(2)}\right) \\
&\boldsymbol{\sigma}_{(23)} = \frac{3F\left(1-\boldsymbol{x}_{(2)}^2\right)}{8} + \frac{F\nu\left(3\boldsymbol{x}_{(1)}^2-1\right)}{8(1+\nu)} \\
&\qquad\qquad - \frac{3F\nu}{2\pi^2(1+\nu)} \sum_{n=1}^{\infty} \frac{(-1)^n}{n^2\cosh(n\pi)} \cos\left(n\pi \boldsymbol{x}_{(1)}\right) \cosh\left(n\pi \boldsymbol{x}_{(2)}\right)
\end{aligned} \tag{91}$$

The displacement fields corresponding to these stresses, up to the addition of a rigid body motion, is given by

$$\begin{aligned}
\boldsymbol{u}_{(1)} &= -\frac{3F\nu}{4E_{\mathsf{Y}}}\boldsymbol{x}_{(1)}\boldsymbol{x}_{(2)}\boldsymbol{x}_{(3)} \\
\boldsymbol{u}_{(2)} &= \frac{F}{8E_{\mathsf{Y}}}\left[3\nu\boldsymbol{x}_{(3)}\left(\boldsymbol{x}_{(1)}^2-\boldsymbol{x}_{(2)}^2\right) - \boldsymbol{x}_{(3)}^3\right] \\
\boldsymbol{u}_{(3)} &= \frac{F}{8E_{\mathsf{Y}}}\left[3\boldsymbol{x}_{(2)}\boldsymbol{x}_{(3)}^2 + \nu\boldsymbol{x}_{(2)}\left(\boldsymbol{x}_{(2)}^2 - 3\boldsymbol{x}_{(1)}^2\right)\right] + \frac{2(1+\nu)}{E_{\mathsf{Y}}}z(\boldsymbol{x})
\end{aligned} \tag{92}$$

where $z(\boldsymbol{x})$ is the anti-derivative of $\boldsymbol{\sigma}_{(23)}$ with respect to $\boldsymbol{x}_{(2)}$. In the present numerical study, the length of the beam is $L = 10$, the shear load is taken as $F = 0.1$, and the material properties are selected as $E_{\mathsf{Y}} = 25$ and $\nu = 0.3$. Moreover, $\Gamma_{\boldsymbol{u}}$ is taken to be the face passing through $\boldsymbol{x}_{(3)} = L$ and boundary displacements are set to $\boldsymbol{g} = \boldsymbol{u}|_{\Gamma_{\boldsymbol{u}}}$. In addition to CVT and RND meshes, we also consider uniform meshes of hexahedral (brick) elements. This also allows for a direct comparison with standard trilinear finite elements on hexahedra. As a means to visualize the beam deformation and the stress field generated under the shear load, we show one set of our results in Fig. 9 for CVT and RND meshes.

In our numerical studies, we have observed that the use of boundary coordinates of [3] with the known surface integrals yields almost identical solutions (less than 0.2% and 0.4% difference in the errors for CVT and RND meshes, respectively) to those obtained from the application of nodal quadrature for computing the projection maps. Subsequently, we will only present results using the nodal quadrature in the remainder of this section.

First, we verify the convergence of the VEM under refinement for different mesh types. The results of this study are shown in Fig. 10 where the relative errors are plotted against the average element diameters. Here the scaling coefficient is set to $\alpha^E = \alpha_\star^E$ corresponding to $\gamma = 1$. For CVT and RND meshes, each point in the curve represents the average values of five sets of meshes with the same number of elements. As evident from the plots, we have second-order convergence in displacements and first-order convergence in stresses in all cases. The fact that these optimal convergence rates are observed for RND meshes with many irregular elements is encouraging and a testament to the robustness of the method.

As a way of comparing the performance of the VEM for the different mesh type, we next plot the error as a function of total number of degrees of freedom (DoFs) in Fig. 11. For the displacement fields, the CVT and RND meshes produce similar errors, while the hexahedral meshes are marginally more accurate. The same trend is observed for the errors in the stress field though the difference between the mesh types is more pronounced.

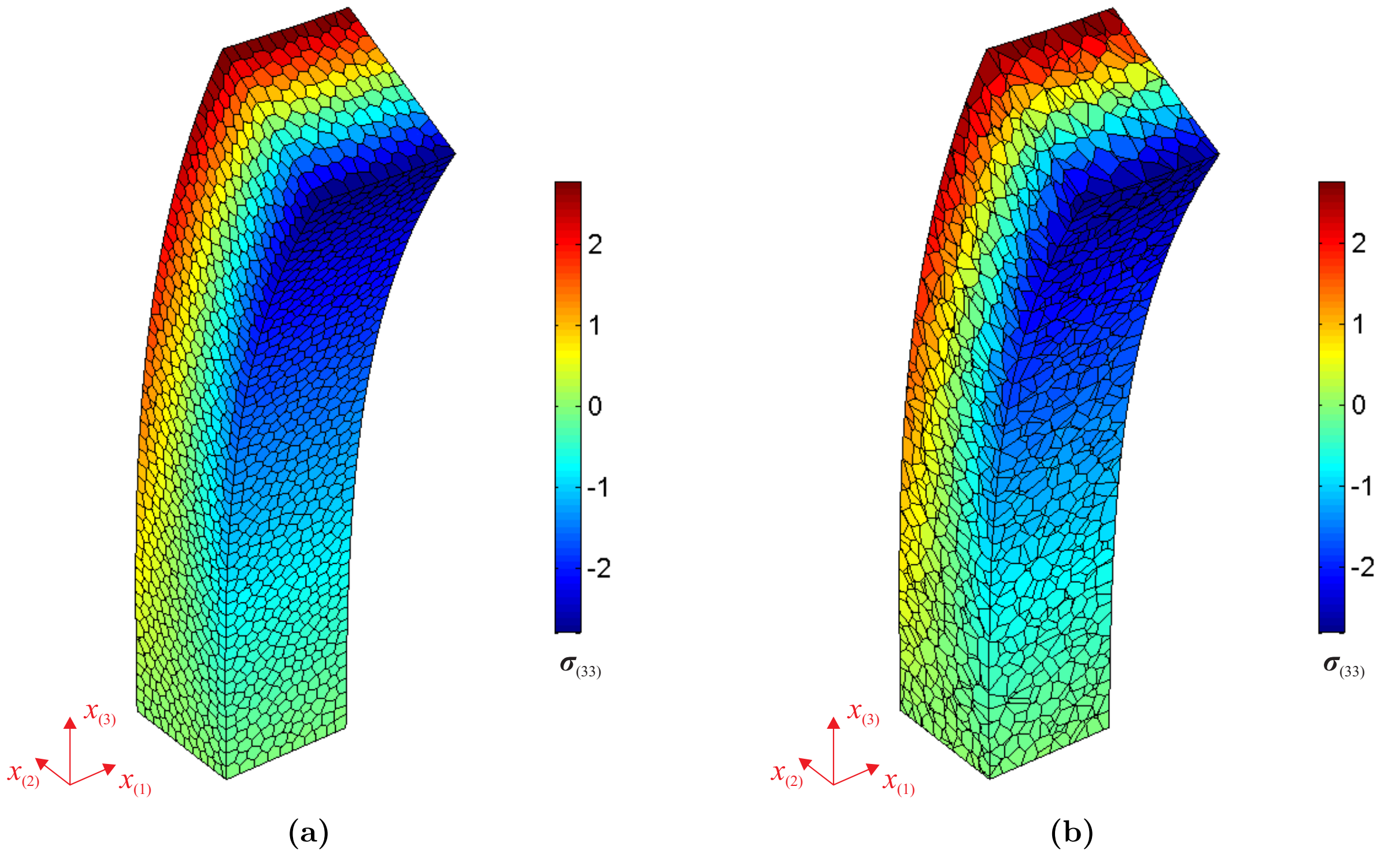

**Figure 9:** Deformation plots for shear loaded cantilever beam bending problem for representative (a) CVT and (b) RND meshes. The colors indicate the magnitude of $[\boldsymbol{\sigma}_h(\tilde{\boldsymbol{u}}_h)]_{(33)}$. The shear load is applied on the bottom face in the negative $\boldsymbol{x}_{(2)}$ direction.

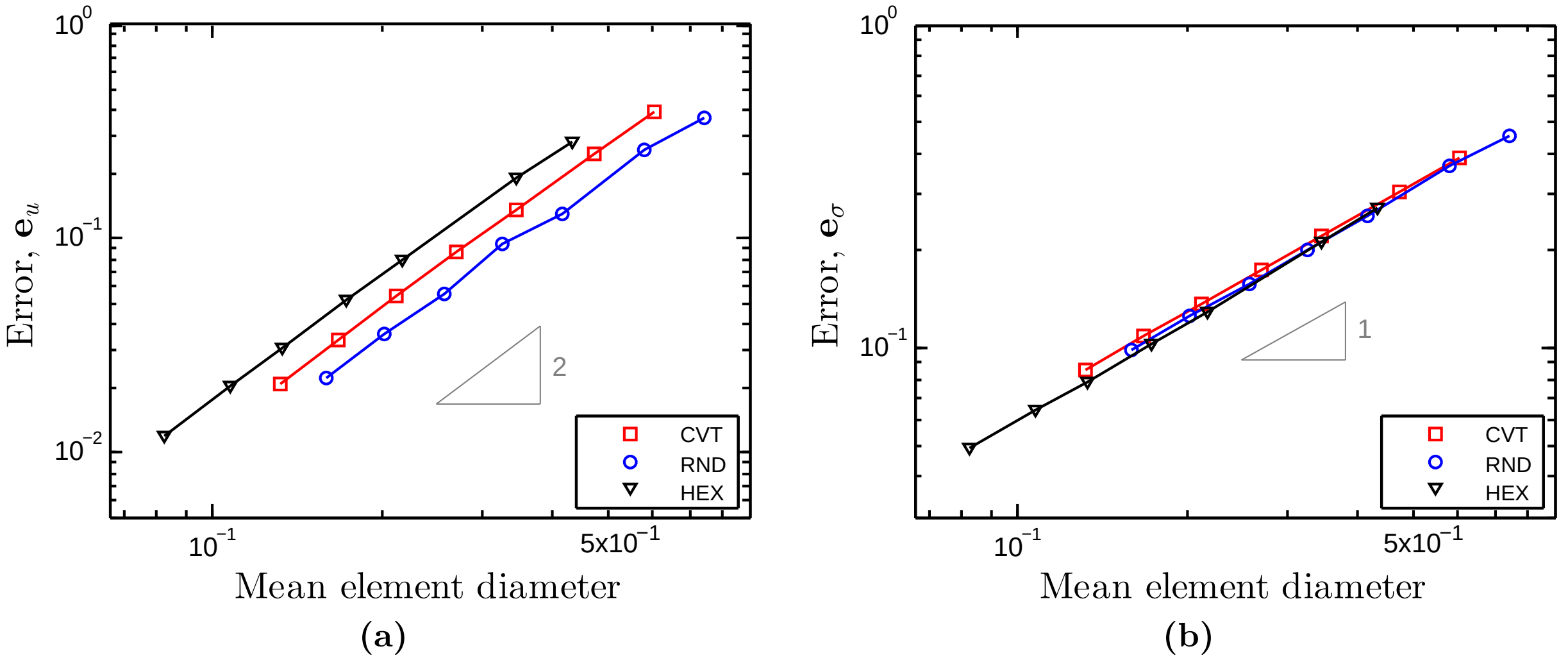

**Figure 10:** Convergence study for VEM under mesh refinement on different meshes, namely CVT Voronoi, random Voronoi (RND) and hexahedral mesh (HEX). (a) Displacement errors, $\mathsf{e}_{\boldsymbol{u}}$. (b) Stress errors, $\mathsf{e}_{\boldsymbol{\sigma}}$. Results pertaining to polyhedral meshes are average of 5 meshes.

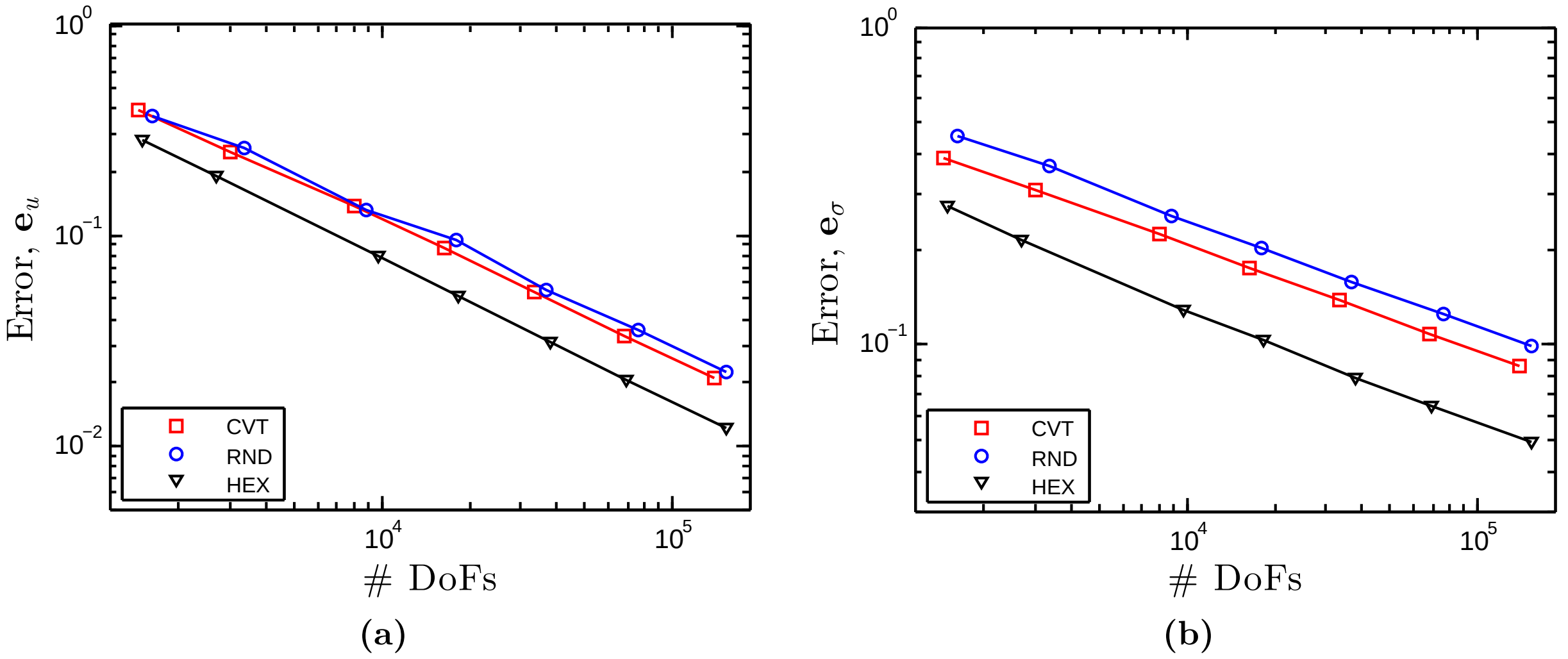

**Figure 11:** Plot of error versus total number of degrees of freedom (DoFs) for CVT, RND, and HEX meshes (a) Displacement errors, $\mathsf{e}_{\boldsymbol{u}}$. (b) Stress errors, $\mathsf{e}_{\boldsymbol{\sigma}}$. Results pertaining to polyhedral meshes are the average of 5 mesh realizations.

While the preceding theoretical discussion illustrates that optimal convergence of the solutions can be expected with any fixed value of $\gamma$ for sufficiently regular meshes, the choice of $\gamma$ can play a significant role in the accuracy of solutions. In order to study this effect, we next perform a parametric study and plot errors for the shear-loaded beam as a function of $\gamma$ for each mesh type (see Fig. 12). The results for CVT and RND meshes are the average of five meshes with 200 elements. The hexahedral mesh has 625 elements but comparable number of degrees of freedom to the polyhedral meshes. As a point of reference, the errors produced by classical trilinear finite elements on this hexahedral mesh are also indicated on the plots.

A few observations regarding these results are in order. First, as seen in Fig. 12(a), there exists an optimal value of $\gamma$ for each mesh type where the displacement errors are minimized. The solutions can be significantly more accurate for this value of $\gamma$. For example, the VEM solution on the hexahedral mesh has almost two orders of magnitude smaller error compared to corresponding finite element solution. We also note that the sensitivity of error is larger for $\gamma$ values less than this optimal value. Conversely, there is less variation in errors for larger values of $\gamma$ across mesh types.

In contrast to the displacement errors, the stress errors are far less sensitive to the change in $\gamma$ and almost constant in the range considered here. This could be partially attributed to the fact that the influence of higher-order modes do not show up in this measure of error. We also note that the VEM error levels on the hexahedral mesh are comparable to the finite element errors.

In the final study, summarized in Fig. 13, we study the effects of mesh refinement on the optimal value of the scaling coefficient. Here the optimal value of $\gamma$, with respect to $\mathsf{e}_{\boldsymbol{u}}$, are determined and plotted for different meshes at different levels of refinement.

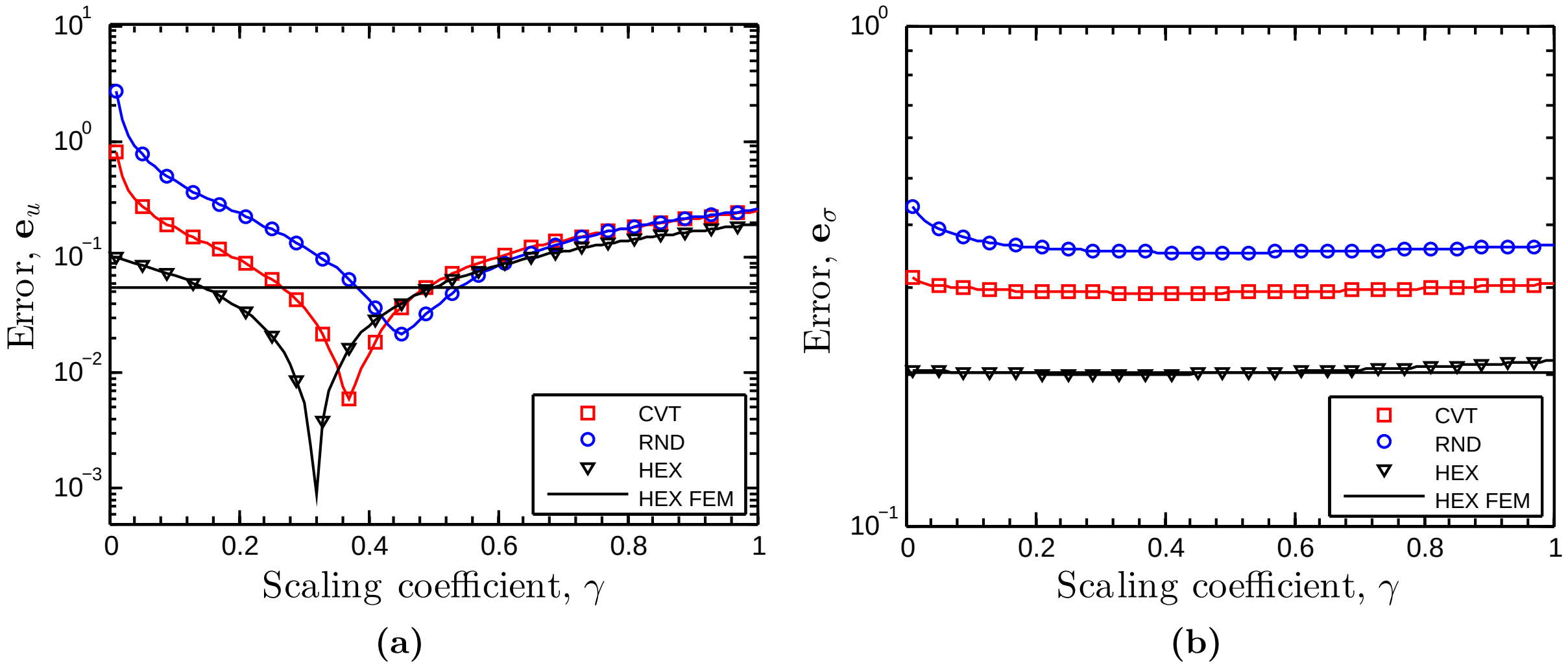

**Figure 12:** Study of the optimal scaling coefficient $\gamma$ for the shear-loaded beam problem for different meshes. (a) Displacement errors, $\mathsf{e}_{\boldsymbol{u}}$. (b) Stress errors, $\mathsf{e}_{\boldsymbol{\sigma}}$. Each point of the CVT and RND curves represents the average errors of 5 meshes.

Aside from a mild increase for the CVT meshes on average, the optimal $\gamma$ is fairly stable under mesh refinement. It may be possible to exploit this fact in practice and use coarse meshes to estimate an optimal range for the scaling coefficient for the problem at hand, subsequently used for a more refined analysis.

We stress, however, that the optimal value of $\gamma$ observed here in these studies are particular to the shear-load beam problem and one cannot make generalizations for other problems. Motivated by element-level energetic considerations of Section 4.4, we recommend using $\gamma = 1$ for general polyhedral meshes in the absence of additional insights into the nature of the problem at hand.

## 6. Concluding remarks

In this work, we discussed the theoretical and practical aspects of VEM in the context of three-dimensional elasticity problems on polyhedral meshes. At the core of the method, local polynomial projections are used to decompose the element strain energy into its uniform strain and higher-order components. When constructing an approximate (discrete) strain energy for an element, preserving the former guarantees the satisfaction of the patch test and the consistency of the method. One can ensure first-order convergence with a suitable but possibly crude approximation of the energy of higher-order modes. As discussed in [37], this splitting of energy can be useful for restoring consistency for polyhedral finite elements in the presence of quadrature error in the evaluation of strain energy. As such, the core concepts underlying VEM can provide an alternative approach for addressing the challenges facing finite element schemes on arbitrary grids.

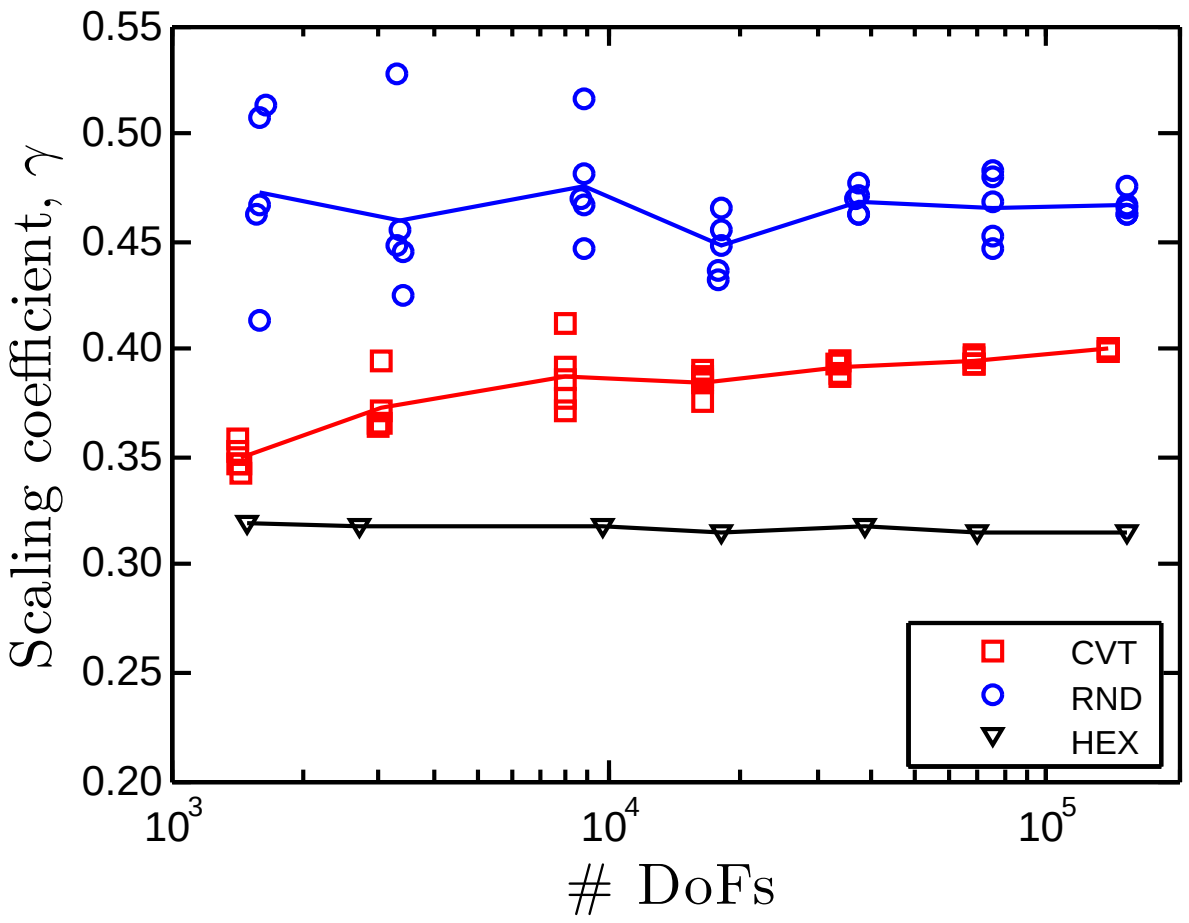

**Figure 13:** Illustration of optimum $\gamma$ for the shear-loaded beam problem under refinement for CVT, RND, and HEX meshes.

## Acknowledgements

We are thankful to the support from the US National Science Foundation under grant number 1321661, and from the Donald B. and Elizabeth M.Willett endowment at the University of Illinois at Urbana-Champaign. We also acknowledge partial support provided by Tecgraf/PUC-Rio (Group of Technology in Computer Graphics), Rio de Janeiro, Brazil. Any opinion, finding, conclusions or recommendations expressed here are those of the authors and do not necessarily reflect the views of the sponsors or sponsoring agencies. We acknowledge Prof. Anton Evgrafov for suggesting, during the World Congress of Computational Mechanics (WCCM 2012, Brazil), that we investigate the class of mimetic numerical methods in order to develop our polygonal discretization scheme in three dimensions – his suggestion was the inception of the work in the present paper.

## Appendix

In the recent work by Ahmed et al. [3], a particular set of barycentric coordinates on polygons are constructed such that their first-order moments can be computed exactly. We will now briefly discuss these coordinates, which will be in turn used in the construction of local element space $\mathcal{W}(E)$.

Consider a polygon $F \subseteq \Re^2$ with $m$ vertices located in counter-clockwise order at $\boldsymbol{x}_1^F, \ldots, \boldsymbol{x}_m^F$. The barycentric coordinate $\varphi_i$ associated with $i$th vertex is equal to one at $\boldsymbol{x}_i$, decays linearly along the incident edges and vanishes at other vertices of $F$. In the interior of the polygon, $\varphi_i$ is a function whose Laplacian is a linear field subject to the condition that its zeroth and first moment equal the corresponding moments of its polynomial projection $p_i$ defined by

$$p_i \doteq \frac{1}{m} + \left( \frac{1}{|F|} \int_{\partial F} \varphi_i \boldsymbol{n} \mathrm{d}s \right) \cdot \left( \boldsymbol{x} - \hat{\boldsymbol{x}}^F \right) \tag{93}$$

Here, $\boldsymbol{n}$ is the unit normal to the boundary of $F$ and $\hat{\boldsymbol{x}}^F = \frac{1}{m} \sum_{j=1}^{m} \boldsymbol{x}_j^F$. Denoting by $\boldsymbol{n}_i$ the unit normal to the edge connecting the $(i-1)$ and $i$th vertices, we can simplify the expression for $p_i$ as

$$p_i = \frac{1}{m} + \frac{1}{2\,|F|}\left(\left|\boldsymbol{x}_i^F - \boldsymbol{x}_{i-1}^F\right|\boldsymbol{n}_i + \left|\boldsymbol{x}_{i+1}^F - \boldsymbol{x}_i^F\right|\boldsymbol{n}_{i+1}\right)\cdot\left(\boldsymbol{x} - \hat{\boldsymbol{x}}^F\right) \tag{94}$$

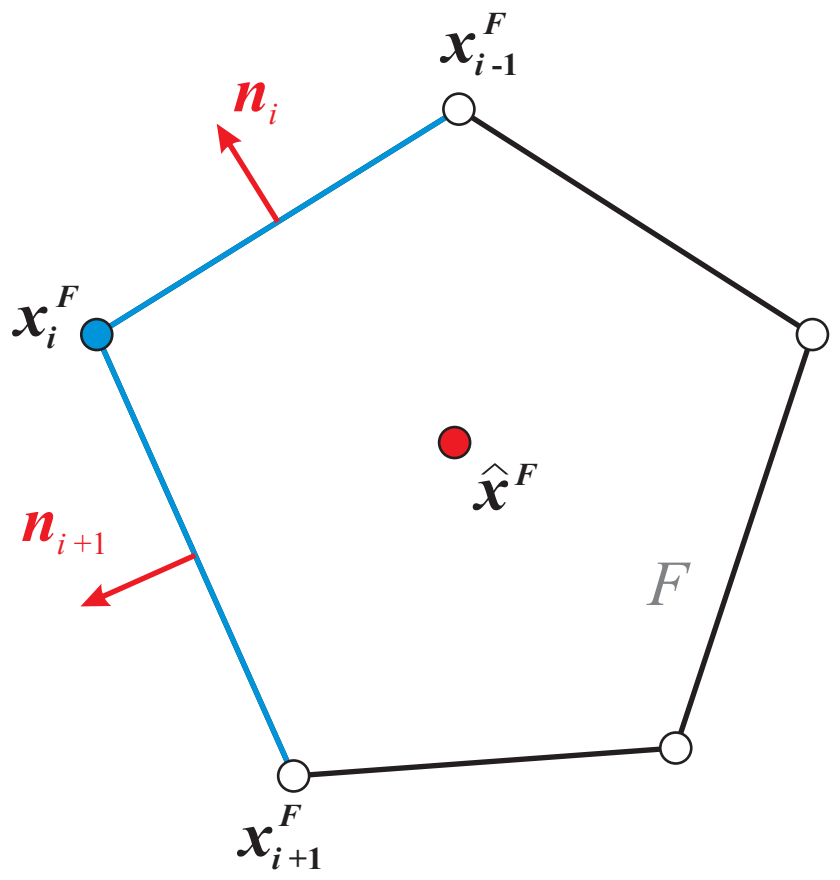

**Figure 14:** Illustration of quantities needed for computing the moment of barycentric coordinates of [3]. Here, $\boldsymbol{n}_i$ and $\boldsymbol{n}_{i+1}$ are the unit normal vectors corresponding to the edges incident on vertex $\boldsymbol{x}_i$ and in the plane of $F$.

Since there is a one-to one mapping between the moments of these functions and their Laplacian, the condition that $\varphi_i$ and $p_i$ have identical first-order moments uniquely defines the coordinates $\varphi_i$. Note that we will not need to compute these coordinates explicitly in the interior of faces in VEM since we only need the average value. This quantity, by construction, can be computed as

$$\int_F \varphi_i \mathrm{d}\boldsymbol{x} = \int_F p_i \mathrm{d}\boldsymbol{x} = \frac{|F|}{m} + \frac{1}{2}\left(\left|\boldsymbol{x}_i^F - \boldsymbol{x}_{i-1}^F\right|\boldsymbol{n}_i + \left|\boldsymbol{x}_{i+1}^F - \boldsymbol{x}_i^F\right|\boldsymbol{n}_{i+1}\right)\cdot\left(\boldsymbol{x}^F - \hat{\boldsymbol{x}}^F\right) \tag{95}$$

where $\boldsymbol{x}^F$ is the centroid of $F$ (cf. Fig. 14).